\documentclass[review]{elsarticle}
\usepackage{placeins}

\usepackage[hyphens]{url}

\usepackage[displaymath]{lineno}
\usepackage{hyperref}

\usepackage{graphicx,pstricks}
\usepackage{graphics}
\usepackage{epsfig}
\usepackage{amssymb}
\usepackage{amsmath, amsthm}
\usepackage{amsfonts}
\usepackage{algorithm}
\usepackage{tabu}
\usepackage{multirow}
\usepackage{subfig}
\usepackage{caption}
\usepackage{changepage}
\usepackage{algpseudocode}

\journal{}

\begin{document}

\begin{frontmatter}

\title{Bayesian Deep Learning for Partial Differential Equation Parameter Discovery with Sparse and Noisy Data}

\author[cee]{Christophe Bonneville\corref{correspondingauthor}}
\author[cee,cam]{Christopher Earls}
\cortext[correspondingauthor]{Corresponding Author: \href{mailto:cpb97@cornell.edu}{\texttt{cpb97@cornell.edu}}}
\address[cee]{School of Civil \& Environmental Engineering, Cornell University, Ithaca, NY 14850, United States}
\address[cam]{Center for Applied Mathematics, Cornell University, Ithaca, NY 14850, United States}

\begin{abstract}
Scientific machine learning has been successfully applied to inverse problems and PDE discovery in computational physics. One caveat concerning current methods is the need for large amounts of (``clean") data, in order to characterize the full system response and discover underlying physical models. Bayesian methods may be particularly promising for overcoming these challenges, as they are naturally less sensitive to the negative effects of sparse and noisy data. In this paper, we propose to use Bayesian neural networks (BNN) in order to: 1) Recover the full system states from measurement data (e.g. temperature, velocity field, etc.). We use Hamiltonian Monte-Carlo to sample the posterior distribution of a deep and dense BNN, and show that it is possible to accurately capture physics of varying complexity, without overfitting. 2) Recover the parameters instantiating the underlying partial differential equation (PDE) governing the physical system. Using the trained BNN, as a surrogate of the system response, we generate datasets of derivatives that are potentially comprising the latent PDE governing the observed system and then perform a sequential threshold Bayesian linear regression (STBLR), between the successive derivatives in space and time, to recover the original PDE parameters. We take advantage of the confidence intervals within the BNN outputs, and introduce the spatial derivatives cumulative variance into the STBLR likelihood, to mitigate the influence of highly uncertain derivative data points; thus allowing for more accurate parameter discovery. We demonstrate our approach on a handful of example, in applied physics and non-linear dynamics.
\end{abstract}


\end{frontmatter}


\section{Introduction}


\noindent In recent years, pioneering research has been conducted into the application of machine learning to computational physics and engineering contexts: example works include \cite{RAISSI2017683, RAISSI2018125, RAISSI2019686, Brunton3932, Rudye1602614, Lusch_2018}. As a result, new sub-fields within the computational sciences have emerged, including, but not limited to physics-informed machine learning \cite{nature_karniadakis} and scientific machine learning (SciML) \cite{osti_1478744}. Within these sub-fields, the development of machine learning based methods to infer the parameters of dynamical system governing equations, and/or the discovery of new partial differential equations (PDE) has attracted significant attention and important early success. In such early work, different inference strategies have been proposed. A popular method, known as \textit{SINDy} \cite{Brunton3932} has its foundations in building a dataset of spatial derivatives that are potentially involved in the governing equation of the observed physics, in order to perform a sparse linear regression aimed at estimating the coefficients of each derivative term forming some latent PDE. This seminal work was later extended in \cite{Rudye1602614}, where the PDE derivative terms were computed either through polynomial interpolations or finite differences, and in \cite{Xu_2021, XU2021110592, Both_2021} using neural networks. A major advantage of these methods is the interpretability: the spatial derivatives involved in the PDE, along with their coefficients, are discovered explicitly. Other methods directly approximate the differential operator using a physics-informed neural network (PINN) \cite{RAISSI2019686}. While these methods generally yield highly accurate forward models, they lack interpretability \cite{JMLR:v19:18-046}, or only allow for recovering the coefficients of a PDE with the functional form already known \cite{RAISSI2019686}. 
\\\\
\noindent The sparse regression based methods outlined earlier have been providing promising results, but the regression system is often poorly conditioned, and the spatial derivatives estimations must be highly accurate in order to provide satisfying results. In \cite{Xu_2021, XU2021110592, Both_2021}, measured data from a physical system are interpolated using a neural network, and then differentiated to create the spatial derivative dataset. This requires abundant data, with as little noise as possible. In practical engineering and scientific applications however, acquiring enough data (which is typically experimental measurements) may be a very expensive and time consuming process. While deep and dense neural networks are able to capture complex physics described in snapshots of physical system dynamics (e.g. shocks or sharp gradients), they are also more likely to overfit noisy measurements. Unfortunately, methods for fine tuning the networks parameters (e.g. cross validation) may be inapplicable due to the lack of abundant data in many scientific applications. Consequently, there is a trade-off: Use more data for training, and potentially obtain better derivative estimations, but with a significant risk of overfitting, or use more data for testing and limit overfitting, but incur a greater risk of missing complex underlying physics. 
\\\\
\noindent Bayesian methods, in particular Bayesian neural network (BNN) \cite{10.5555/525544, 10.5555/1162264, Goan_2020}, are promising for avoiding this trade-off; as they a naturally less prone to overfitting, even with very noisy and sparse data. Furthermore, BNNs provide confidence intervals on their prediction, which can be used for improving the accuracy of PDE discovery techniques. In the field of PDE discovery, Bayesian machine learning methods have been used to infer parameters of governing equations with known functional form. For example, \cite{RAISSI2017683} used Gaussian processes, and \cite{Yang_2021, MENG2021110361} used BNNs combined with PINNs to infer PDE parameters, but do not permit the discovery of unknown PDEs. In this paper, we extend the original approach outlined in \cite{Brunton3932, Rudye1602614, Xu_2021, Both_2021, zhang_2018} in two important ways. Firstly, we propose to use a BNN to interpolate the snapshot measurements from the physical system. Relying on BNNs offers two major advantages: It makes the network more robust to noise and sparse data (allowing for minimal tuning of the network's hyperparameters) and it provides valuable confidence intervals over the interpolation predictions. Secondly, we propose to use these BNN confidence intervals to quantify the uncertainty over the spatial derivative dataset, and introduce this uncertainty into a sequential threshold Bayesian linear regression model (STBLR) to recover the PDE coefficients.
\\\\
\noindent In the following sections, we first provide elements of background on neural networks and Bayesian inference (section \ref{section_bnn}). We then introduce our framework for PDE discovery (section \ref{section_framework}) and present our results on three application examples (section \ref{section_example}).

\section{Bayesian Neural Network Background}
\label{section_bnn}

\subsection{Standard Deep Neural Networks}

\begin{figure}[!h]
\centering
  \centering
  \includegraphics[width=\textwidth]{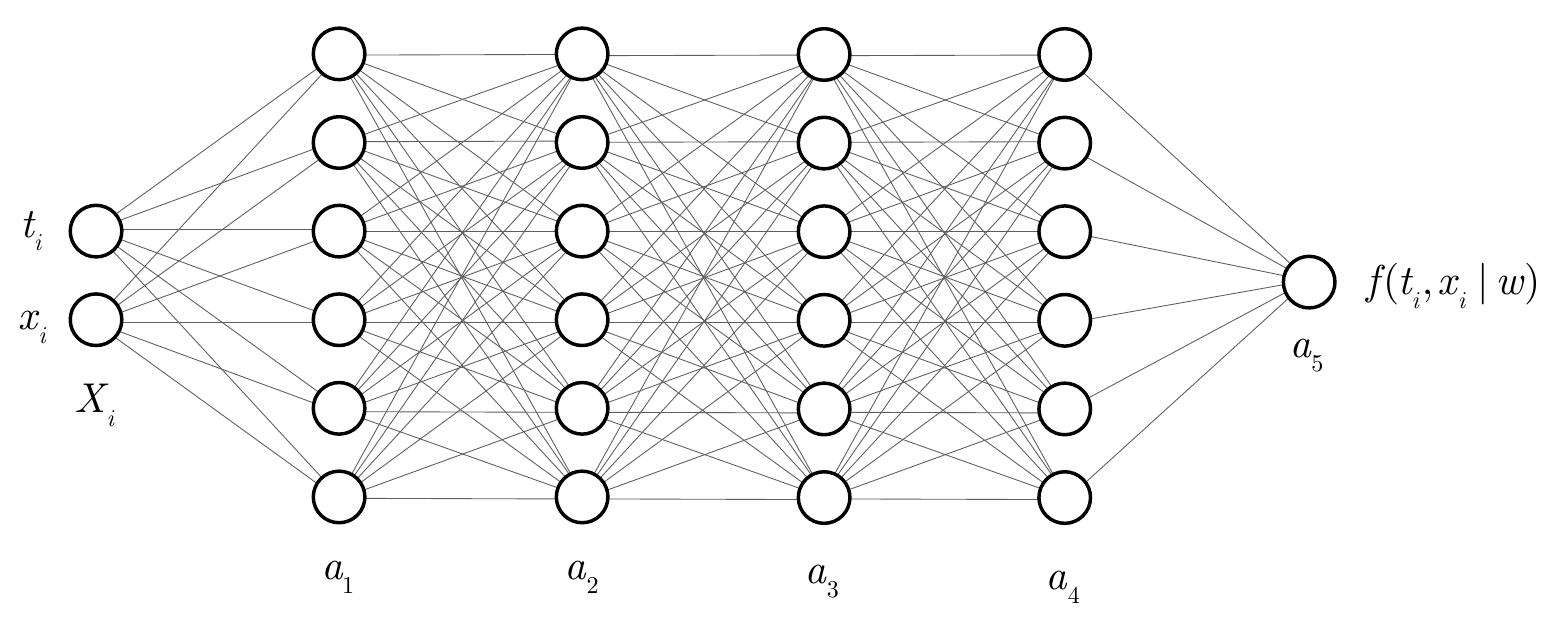}
  \caption{Representation of a neural network with 4 hidden layers and 6 hidden units}
  \label{fig_background1}
\end{figure}

\noindent A neural network \cite{10.5555/1162264, murphy2013machine, GoodBengCour16} is a non-linear parametric function, able to learn and approximate any other continuous function representations under weak conditions, provided that the network is sufficiently complex \cite{HornikEtAl89}. Figure \ref{fig_background1} shows a standard fully connected neural network architecture with $n_l=4$ hidden layers and $n_u=6$ units, respectively. Analytically, the output of a neural network for regression with a time-space input coordinate $X_i=(t_i,x_i)$ can be written as a function composition:
\begin{equation}
\begin{aligned}
    f(t_i,x_i|w)&=a_5=w_5\cdot a_4=w_5\cdot\varphi(w_4\cdot a_3)\\
    &=w_5\cdot\varphi(w_4\cdot (\ldots \cdot \varphi(w_1\cdot X_i))
\end{aligned}
\end{equation}

\noindent Where $\varphi$ is a non-linear activation function, and the set of parameters (weights) is written $w=\{w_l\}$ with $l\in[\![1,n_l+1]\!]$. In a neural network regression framework, we may further assume that the known data have been generated from the network, itself, and subsequently corrupted with Gaussian noise \cite{10.5555/1162264, murphy2013machine}. That is, for input $X_i=(t_i,x_i)$, we have:
\begin{equation}
y_i=f(t_i,x_i|w)+\epsilon\hspace{1cm}\epsilon\sim\mathcal{N}(0,\sigma^2)
\end{equation}
Leading to the well-known likelihood function:
\begin{equation}
\label{likelihood}
    p(y|X,w)=\mathcal{N}(y|f(X|w),\sigma^2I)
\end{equation}
\noindent In standard non-Bayesian deep learning applications, we are generally interested in maximizing equation \ref{likelihood} (or variations of it) with respect to the weights $w$, using numerical optimization methods \cite{GoodBengCour16}.

\subsection{Bayesian Inference for Neural Networks}

\noindent In a BNN, the weights are assumed to be sampled from a prior probability distribution, for example a standard Gaussian distribution with a 0 mean and unit standard deviation:
\begin{equation}
\label{prior}
    w\sim\mathcal{N}(0, I)
\end{equation}
Using equation \ref{likelihood} and Bayes' rule, the posterior distribution over the weights can be written as:
\begin{equation}
\begin{aligned}
    p(w|X,y)&=\frac{p(y|X,w)p(w)}{p(y|X)}\\[5pt]
    &\propto\mathcal{N}(y|f(X|w),\sigma^2I)\mathcal{N}(w|0, I)
\end{aligned}
\end{equation}
Due to the highly non-linear nature of $f(t,x|w)$, the likelihood is a very complicated function; thus sampling from the posterior is analytically intractable. We can either approximate the posterior with e.g. a Gaussian distribution and find the parameters of the Gaussian by minimizing the KL divergence between the true and the approximate posterior (variational inference) \cite{NIPS2011_7eb3c8be, Goan_2020, jospin2020handson, Blei_2017}, or use Monte-Carlo sampling \cite{10.5555/525544, cobb2020scaling, jospin2020handson, Goan_2020}. Variational inference for BNNs generally offers poor accuracy and often fails to provide meaningful confidence intervals \cite{Goan_2020, Blei_2017, pmlr-vR5-wang05a, turner_sahani_2011}. Here, we rely on Hamiltonian Monte-Carlo (HMC) instead \cite{10.5555/525544, betancourt2018conceptual, Neal2011}\\\\
\noindent HMC flips the log-posterior probability distribution to be ``up-side-down", so that the regions of high probability become minima. If a virtual particle (with momentum, $v$) is placed on the flipped distribution and starts rolling freely, it will naturally go towards the regions of lower potential energy (minima), i.e. regions with higher probability. HMC has two steps. First it simulates the motion of such a particle, using Hamiltonian physics, and subsequently records the position over time; thus providing a set of sample candidates. The method then uses Metropolis Hastings \cite{doi:10.1063/1.1699114} to either accept or reject these samples. The Hamiltonian step generally provides suitable samples, and allows for a higher acceptance rate than other Markov-Chain-Monte-Carlo (MCMC) methods \cite{10.5555/525544}. The Hamiltonian can be expressed as a joint probability between $w$ and $v$:
\begin{equation}
\begin{aligned}
    \mathcal{H}(v,w)&=-\log(p(v,w|X,y))\\
    &=-\log(p(v|w,X,y))-\log(p(w|X,y))\\
    &=T(v|w)+V(w)
\end{aligned}
\end{equation}
Where $T$ and $V$ are the kinetic and potential energy, respectively. Using Hamilton's equations, and assuming that the momentum is independent from $w$, this leads to:
\begin{equation}
\label{hamil}
    \frac{dw}{d\tau}=\frac{\partial \mathcal{H}}{\partial v}=\frac{\partial T}{\partial v}\hspace{1cm}\frac{dv}{d\tau}=-\frac{\partial \mathcal{H}}{\partial w}=-\frac{\partial V}{\partial w}
\end{equation}
Equation \ref{hamil} can then be integrated using standard integration algorithms (Euler, Leap-Frog, etc. \cite{Neal2011}) to find sample candidates.

\section{Inversion Framework for PDE Discovery}
\label{section_framework}

\noindent We consider $n$ measurements from a given physical system, whose response is of the form of inputs-outputs, $\mathcal{D}=\{X_i,y_i\}_{i\in[\![1,n]\!]}$; where the input is a time-space coordinate, and the output is a corresponding measured physical quantity (potentially noisy). In this work, all the response data are generated from numerical solutions of PDEs, using Chebyshev polynomials with a Runge-Kutta time integration scheme, using the \textit{Chebfun} matlab package \cite{Driscoll2014}.\\

\noindent Unless specified otherwise, we assume 16 measurement sensors placed at random space locations within the problem domain, that record the physical quantities of interest over time with increment $\Delta t$. The input data have the form $X_i=(t_i,x_i)$, where $t_i$ is time and $x_i$ is the sensor Cartesian coordinate. The output has the form $y_i=u(t_i,x_i)+\epsilon$, where $u$ is the PDE solution, and $\epsilon$ is noise corrupting the data. We consider three cases: $\epsilon=0$ (noiseless), $\epsilon\sim\mathcal{N}(0,0.01^2)$ and $\epsilon\sim\mathcal{N}(0,0.05^2)$.\\

\begin{figure}[!h]
\hspace*{-3.1cm}
  \includegraphics[width=1.5\textwidth]{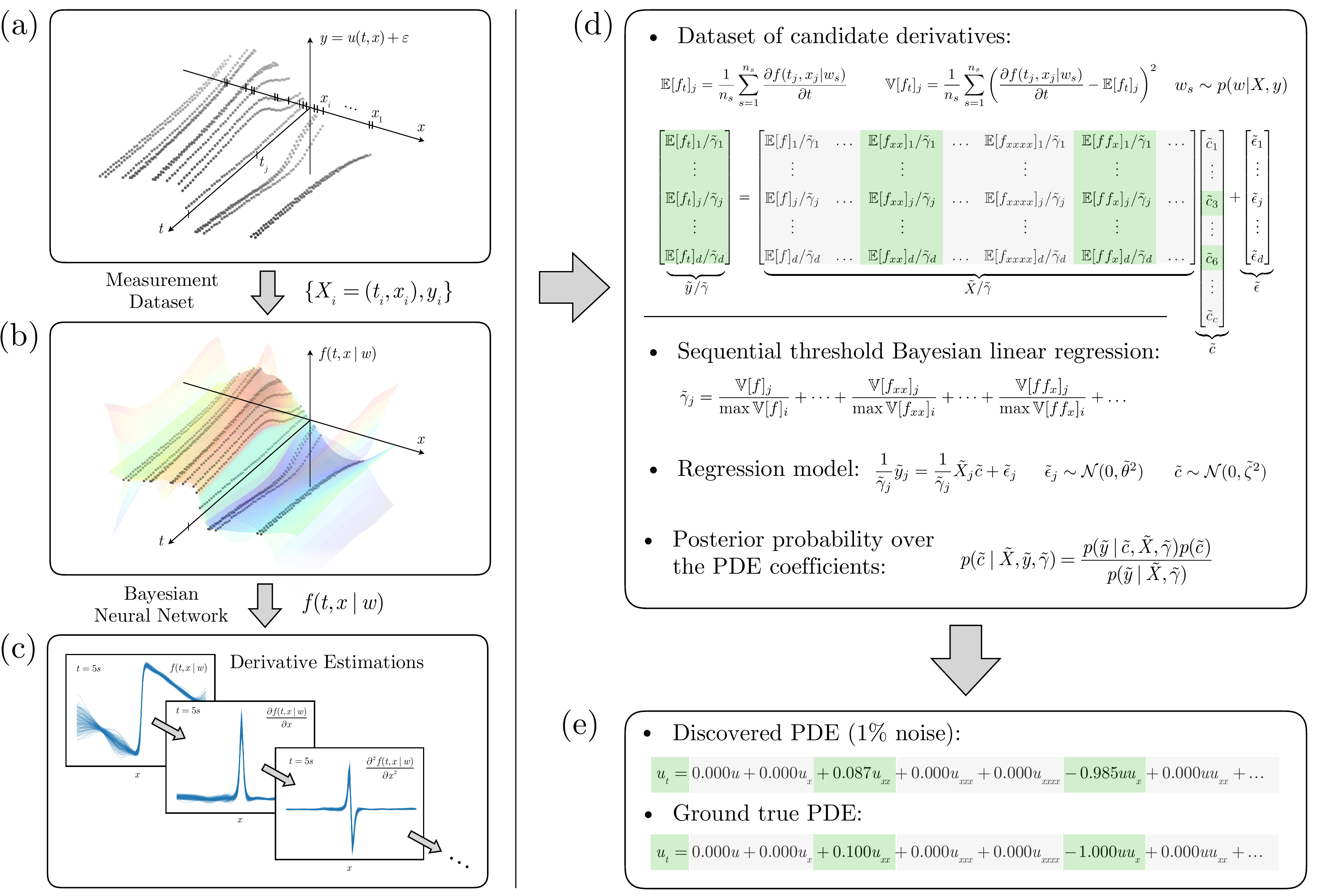}
  \caption{General PDE discovery framework. (a) Measurement data representing noisy snapshots of physical system dynamics, $y_i=u(t_i,x_i)+\epsilon_i$ in time and space. (b) Fitting of the measurement data with a BNN. (c) Differentiation of the trained BNN with respect to time and space, for each sets of weights sampled from the posterior. (d) Construction of a dataset of derivatives, potentially comprising the underlying PDE governing the system, obtained from the BNN. The derivative values are stochastic, and we take their expected value with respect to the BNN posterior. A sequential threshold Bayesian linear regression is performed on the derivative dataset, weighted by the derivative variance, to obtain the value of each PDE coefficients. (e) Discovery of the coefficients and derivatives involved in the underlying PDE (here Burgers equation with noisy measurement data).}
  \label{fig_framework}
\end{figure}

\noindent The general framework for PDE discovery using these system response data is outlined in the following subsections, and summarized in figure \ref{fig_framework}.

\subsection{Bayesian Neural Network Fitting}
\noindent  We assume that $u$ can be approximated with a Bayesian neural network, $f$ (equation \ref{approx}), trained on the measurement data using HMC.
\begin{equation}
\label{approx}
\begin{aligned}
    y_i&=u(t_i,x_i)+\epsilon\hspace{1cm}\epsilon\sim\mathcal{N}(0,\sigma^2)\\
    &\hspace{1.2cm}u(t,x)\approx f(t,x|w)
\end{aligned}
\end{equation}
In the following sections, the BNN indifferently employs a single architecture: $n_l=4$ fully connected hidden layers of $n_u=50$  hidden units each, with $\varphi=\tanh$ activation functions. Unless specified otherwise, the prior over each neural network weight is a standard Gaussian distribution (zero mean and unit standard deviation). The noise standard deviation in the likelihood function is assumed as $\sigma=0.01$. In the HMC sampling we take $n_s=6000$ samples ($200$ burn), with an integration time stepping $\Delta \tau=5\cdot10^{-4}$ ($\Delta \tau=1\cdot10^{-4}$ in section \ref{kdv}).
\subsection{Surrogate Derivative Dataset}
\noindent Once the BNN is properly trained, we sample a set of $d$ time-space coordinates randomly: $(t_j,x_j)\sim \mathcal{U}[\Omega]$ where $\Omega$ is the time-space domain and $\mathcal{U}$ the uniform distribution (in the following sections, $d=10000$). The BNN is differentiated multiple times with respect to space and time and evaluated at each input $(t_j,x_j)$ with $j\in[\![1,d]\!]$. The derivatives are then averaged over the set of network weights. For example, the expected value and variance of the first order time derivative is:
    
\begin{equation}
\begin{cases}
    \,\,\,\displaystyle\mathbb{E}[f_t]_j=\frac{1}{n_s}\sum_{s=1}^{n_s}\frac{\partial f(t_j,x_j|w_s)}{\partial t}\\[15pt]
    \,\,\,\displaystyle\mathbb{V}[f_t]_j=\frac{1}{n_s}\sum_{s=1}^{n_s}\bigg(\frac{\partial f(t_j,x_j|w_s)}{\partial t}-\mathbb{E}[f_t]_j\bigg)^2
\end{cases}
    \hspace{1cm}w_s\sim p(w|X,y)
\end{equation}

\noindent Using this method and following the approach introduced by \cite{Brunton3932, Rudye1602614}, we build a library of successive expected derivatives in space $\tilde{X}$; potentially comprising the underlying governing PDE. The spatial derivative orders included in $\tilde{X}$ are arbitrary, and we may include non-linear terms as well. The total number of derivative candidates is defined as $n_c$, and in the following section we have $n_c=11$ (the specific list of derivatives used is outlined in table \ref{table2}). Similarly, we compute the corresponding expected time derivative output vector, $\tilde{y}$, and the variance of each spatial derivative (matrix $\tilde{Z}$).
    
\begin{equation}
\Tilde{X}=
    \begin{bmatrix}
\mathbb{E}[f]_1 & \mathbb{E}[f_x]_1 & \dots & \mathbb{E}[f_{xxxx}]_1 & \mathbb{E}[ff_{x}]_1 & \mathbb{E}[ff_{xx}]_1 & \dots\,\,\, \\
\vdots & \vdots & & \vdots & \vdots & \vdots\\
\mathbb{E}[f]_d & \mathbb{E}[f_x]_d & \dots & \mathbb{E}[f_{xxxx}]_d & \mathbb{E}[ff_{x}]_d & \mathbb{E}[ff_{xx}]_d & \dots\,\,\, \end{bmatrix}_{[d\times n_c]}
\end{equation}
\begin{equation}
\Tilde{Z}=
    \begin{bmatrix}
\mathbb{V}[f]_1 & \mathbb{V}[f_x]_1 & \dots & \mathbb{V}[f_{xxxx}]_1 & \mathbb{V}[ff_{x}]_1 & \mathbb{V}[ff_{xx}]_1 & \dots\,\,\, \\
\vdots & \vdots & & \vdots & \vdots & \vdots\\
\mathbb{V}[f]_d & \mathbb{V}[f_x]_d & \dots & \mathbb{V}[f_{xxxx}]_d & \mathbb{V}[ff_{x}]_d & \mathbb{V}[ff_{xx}]_d & \dots\,\,\, \end{bmatrix}_{[d\times n_c]}
\end{equation}
\begin{equation}
\Tilde{y}=
    \begin{bmatrix}
    \mathbb{E}[f_t]_1\\
    \vdots\\
    \mathbb{E}[f_t]_d
    \end{bmatrix}_{[d\times 1]}
\end{equation}
\\\\
\noindent Higher order differentiation tends to exhibit higher variance. To ensure that the uncertainty of each datapoint is not made overly unbalanced by the higher order terms, each term in the variance matrix $\tilde{Z}$ is normalized, by dividing each columns with their maximum value. Then, every normalized term within each row is summed, thus providing a vector, $\tilde{\gamma}(\tilde{Z})$, quantifying the uncertainty of each derivative data point:
\\\\
\makebox[\textwidth]{\parbox{1.1\textwidth}{
\begin{equation}
\tilde{\gamma}(\tilde{Z})=
\begin{bmatrix}
\tilde{\gamma}_1\\
\vdots\\
\tilde{\gamma}_d
\end{bmatrix}_{[d\times 1]}
=
\begin{bmatrix}
\sum_{j=1}^{n_c}\big(\tilde{Z}_{1j}/\max(\tilde{Z}_{ij})_{i\in[\![1,d]\!]}\big)\\
\vdots\\
\sum_{j=1}^{n_c}\big(\tilde{Z}_{dj}/\max(\tilde{Z}_{ij})_{i\in[\![1,d]\!]}\big)
\end{bmatrix}_{[d\times 1]}
\end{equation}}}

\subsection{Sequential Threshold Bayesian Linear Regression Fitting}
\noindent The library of derivative terms is then fitted with a sequential threshold Bayesian linear regression (STBLR) \cite{zhang_2018}, in order to approximate the coefficients associated with each derivative. This is a sparse regression problem, and it is expected that the set of actual derivatives involved in the ground truth PDE should be much smaller than the entire set of candidate derivatives. Hence most coefficients should be 0. Therefore, a STBLR with Gaussian priors should be well suited, as it will naturally skew most coefficients towards 0; it also provides confidence intervals over the recovered PDE coefficients. Using a STBLR was first proposed in \cite{zhang_2018}, but here we go further, by introducing the uncertainty over each derivative term within the library (quantified through $\tilde{\gamma}(\tilde{Z})$) directly into the BLR model, in order to automatically minimize the influence of highly uncertain derivatives. The regression model can be written as follows:
\\
\begin{equation}
\label{blr_model}
    \frac{1}{\tilde{\gamma}_j}\tilde{y}_j=\frac{1}{\tilde{\gamma}_j}\tilde{X}_j\tilde{c}+\tilde{\epsilon}_j\hspace{1cm}
    \begin{cases}
        \tilde{\epsilon}_j\sim\mathcal{N}(0,\tilde{\theta}^2)\\[2pt]
        \tilde{c}\sim\mathcal{N}(0,\tilde{\zeta}^2)
    \end{cases}
\end{equation}
\\
\noindent The system of equations written explicitly is:
\\\\
\makebox[\textwidth]{\parbox{1.2\textwidth}{
\begin{equation}
\begin{bmatrix}
\mathbb{E}[f_t]_1/\tilde{\gamma}_1\\
\vdots\\
\mathbb{E}[f_t]_d/\tilde{\gamma}_d
\end{bmatrix}=
\begin{bmatrix}
\mathbb{E}[f]_1/\tilde{\gamma}_1 & \dots & \mathbb{E}[f_{xxxx}]_1/\tilde{\gamma}_1 & \dots & \mathbb{E}[ff_{xx}]_1/\tilde{\gamma}_1 & \dots\,\,\, \\
\vdots & & \vdots & & \vdots\\
\mathbb{E}[f]_d/\tilde{\gamma}_d & \dots & \mathbb{E}[f_{xxxx}]_d/\tilde{\gamma}_d & \dots & \mathbb{E}[ff_{xx}]_d/\tilde{\gamma}_d & \dots\,\,\, \\ \end{bmatrix}
\begin{bmatrix}
\tilde{c}_1\\
\vdots\\
\tilde{c}_c
\end{bmatrix}
+
\begin{bmatrix}
\tilde{\epsilon}_1\\
\vdots\\
\tilde{\epsilon}_d
\end{bmatrix}
\end{equation}
}}

\vspace{1cm}
\noindent Each input-output pair $\{\tilde{X}_j,\tilde{y}_j\}$ with $j\in[\![1,d]\!]$ is scaled by $1/\tilde{\gamma}_j$; thus allowing for discounting the influence of highly uncertain derivative data points. Equation \ref{blr_model} leads to the posterior distribution over the PDE coefficients:

\begin{equation}
p(\Tilde{c}|\Tilde{X}, \Tilde{y},\Tilde{\gamma}(\Tilde{Z}))= \frac{p(\Tilde{y}|\Tilde{c},\Tilde{X},\Tilde{\gamma}(\Tilde{Z}))p(\Tilde{c})}{p(\tilde{y}|\tilde{X},\Tilde{\gamma}(\Tilde{Z}))}=\frac{\mathcal{N}(\Tilde{y}/\Tilde{\gamma}|\Tilde{X}\Tilde{c}/\Tilde{\gamma},\Tilde{\theta}^2)\mathcal{N}(\Tilde{c}|0,\Tilde{\zeta}^2)}{\displaystyle \int \mathcal{N}(\Tilde{y}/\Tilde{\gamma}|\Tilde{X}\Tilde{c}/\Tilde{\gamma},\Tilde{\theta}^2)\mathcal{N}(\Tilde{c}|0,\Tilde{\zeta}^2) \text{d}\Tilde{c}}
\end{equation}

\noindent Where $\Tilde{\zeta}$ and $\Tilde{\theta}$ are hyperparameters that maximize the marginal likelihood $p(\tilde{y}|\tilde{X},\Tilde{\gamma}(\Tilde{Z}))$. The vector of coefficients associated with each derivative can be taken as the mean value $\mathbb{E}[\tilde{c}]$ with respect to the posterior (with variance $\mathbb{V}[\tilde{c}]$). Here, Bayesian inference is analytically tractable, and the expected values of interest can be computed exactly.
\\\\
In a STBLR, regression is repeated iteratively. Initially, a wide range of candidate derivatives are assumed so that an initial Bayesian linear regression may be performed. Coefficients having their absolute expected value fall under an arbitrary threshold, $\delta$, are assumed to be, in reality, null; thus the corresponding derivative candidates are removed from the derivative dataset. The process is repeated until all the remaining derivative candidates have absolute expected coefficient values greater than $\delta$. Note that here, we use a dynamic threshold: $\delta$ is doubled at each iteration. STBLR is further detailed in algorithm \ref{alg:cap}.

\begin{algorithm}
\caption{Sequential Threshold Bayesian Linear Regression}\label{alg:cap}
\begin{algorithmic}[1]
\Require $\tilde{X}$, $\tilde{y}$, $\tilde{Z}$, $\delta$
\State $\tilde{X}^{(0)}=\tilde{X}$, $\tilde{Z}^{(0)}=\tilde{Z}$, $\tilde{\gamma}^{(0)}=\tilde{\gamma}(\tilde{Z}^{(0)})$, $\delta^{(0)}=\delta$
\State $\mathbb{E}[\tilde{c}^{(0)}] \gets \text{BLR}(\tilde{X}^{(0)}, \tilde{\gamma}^{(0)}, \tilde{y})$
\State \textbf{find} set of indices $R^{(0)}=\{1,\dots,n_r\}$ for which $|\mathbb{E}[\tilde{c}^{(0)}_j]|<\delta^{(0)}$, $j\in\{1,\dots,n_r\}$
\State $\tilde{X}^{(1)}=\tilde{X}^{(0)}$, $\tilde{Z}^{(1)}=\tilde{Z}^{(0)}$
\For{$j\in R^{(0)}$}
    \State Remove $j^\text{th}$ column in $\tilde{X}^{(1)}$, $\tilde{Z}^{(1)}$
\EndFor
\State $\tilde{\gamma}^{(1)}=\tilde{\gamma}(\tilde{Z}^{(1)})$, $k=1$
\While{$R^{(k-1)}$ is not empty}
    \State $\delta^{(k)}=2\delta^{(k-1)}$
    \State $\mathbb{E}[\tilde{c}^{(k)}] \gets \text{BLR}(\tilde{X}^{(k)}, \tilde{\gamma}^{(k)}, \tilde{y})$
    \State \textbf{find} set of indices $R^{(k)}$ for which $|\mathbb{E}[\tilde{c}^{(k)}_j]|<\delta^{(k)}$, $j\in R^{(k)}$
    \State $\tilde{X}^{(k+1)}=\tilde{X}^{(k)}$, $\tilde{Z}^{(k+1)}=\tilde{Z}^{(k)}$
    \For{$j\in R^{(k)}$}
        \State Remove $j^\text{th}$ column in $\tilde{X}^{(k+1)}$, $\tilde{Z}^{(k+1)}$
        \State $\tilde{\gamma}^{(k+1)}=\tilde{\gamma}(\tilde{Z}^{(k+1)})$, $k=k+1$
    \EndFor
\EndWhile
\end{algorithmic}
\end{algorithm}

\subsection{Error Metric and Computer Implementation}

\noindent For assessing the accuracy of our models, we consider two error metrics. First, the $\ell^2$ norm of the difference between the ground truth and the discovered vector of PDE coefficients:
\begin{equation}
    e_C=|\!|\mathbb{E}[c]-c_\text{true}|\!|_2
\end{equation}
Secondly, we consider the $\ell^2$ norm of the difference between the ground truth PDE solution and the numerical solution of the discovered PDE (we solve it using the \textit{chebfun} package \cite{Driscoll2014}):
\begin{equation}
    e_L=|\!|u_\text{true}(t,x)-u_\text{discovered}(t,x)|\!|_2
\end{equation}

\noindent For baseline comparison, a standard deep neural network (DNN) is also trained with the same architecture (i.e. same hyperparameters) as the BNN, using a mean-squared error loss with a learning rate, $\alpha=2\cdot 10^{-4}$, over $30000$ iterations (we use Adam for gradient descent \cite{DBLP:journals/corr/KingmaB14}). We also compare STBLR with a sequential threshold ordinary least squared regression (STOLS). STOLS follows the same sequential idea as STBLR, but the Bayesian linear regression is replaced with ordinary least squares, and each term in $\tilde{\gamma}(\tilde{Z})$ carries the same weight ($\forall j$, $\tilde{\gamma}_j=1$). 
\\\\
\noindent To implement our neural network models, we use \texttt{PyTorch} \cite{NEURIPS2019_9015}. HMC sampling is performed using the \texttt{Hamiltorch} add-on library \cite{cobb2020scaling}, and the linear regressions for discovering PDE coefficients are performed using \texttt{scikit-learn} \cite{scikit-learn}. Our code and data is available at \href{https://github.com/CBonneville45/BayesianDeepLearningPDE/}{\texttt{github.com/CBonneville45/BayesianDeepLearningPDE}}.
\newpage
\section{Results and Discussion}
\label{section_example}

\subsection{\label{burgers} Burgers Equation}

\begin{figure}[!h]
\centering
  \centering
  \includegraphics[width=1\textwidth]{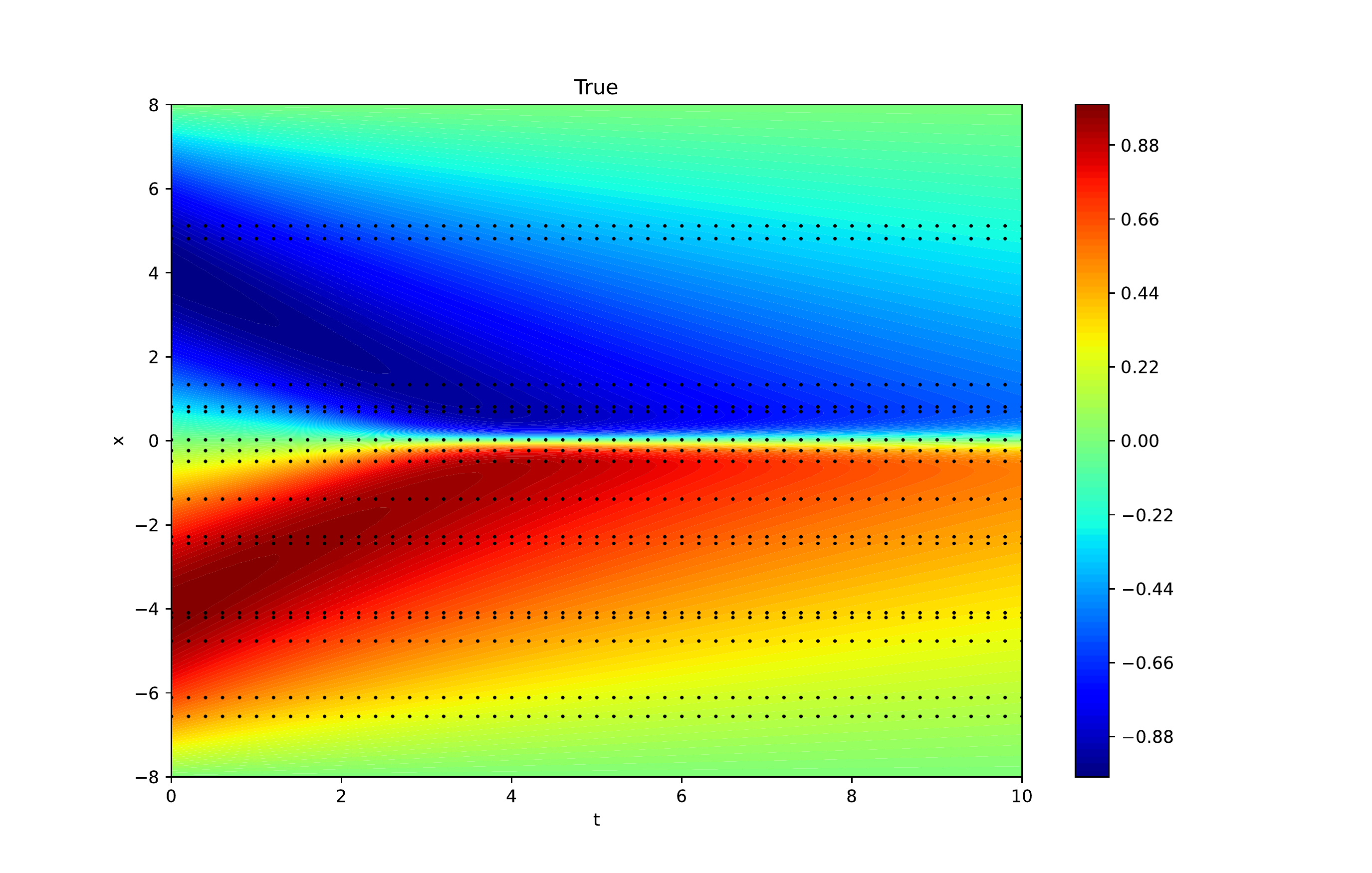}
  \caption{Ground truth solution for Burgers equation (black dots represent measured time-space points)}
  \label{fig4}
\end{figure}

\begin{figure}[!h]
\vspace{-2cm}
\hspace{-2.8cm}
  \includegraphics[width=1.4\textwidth]{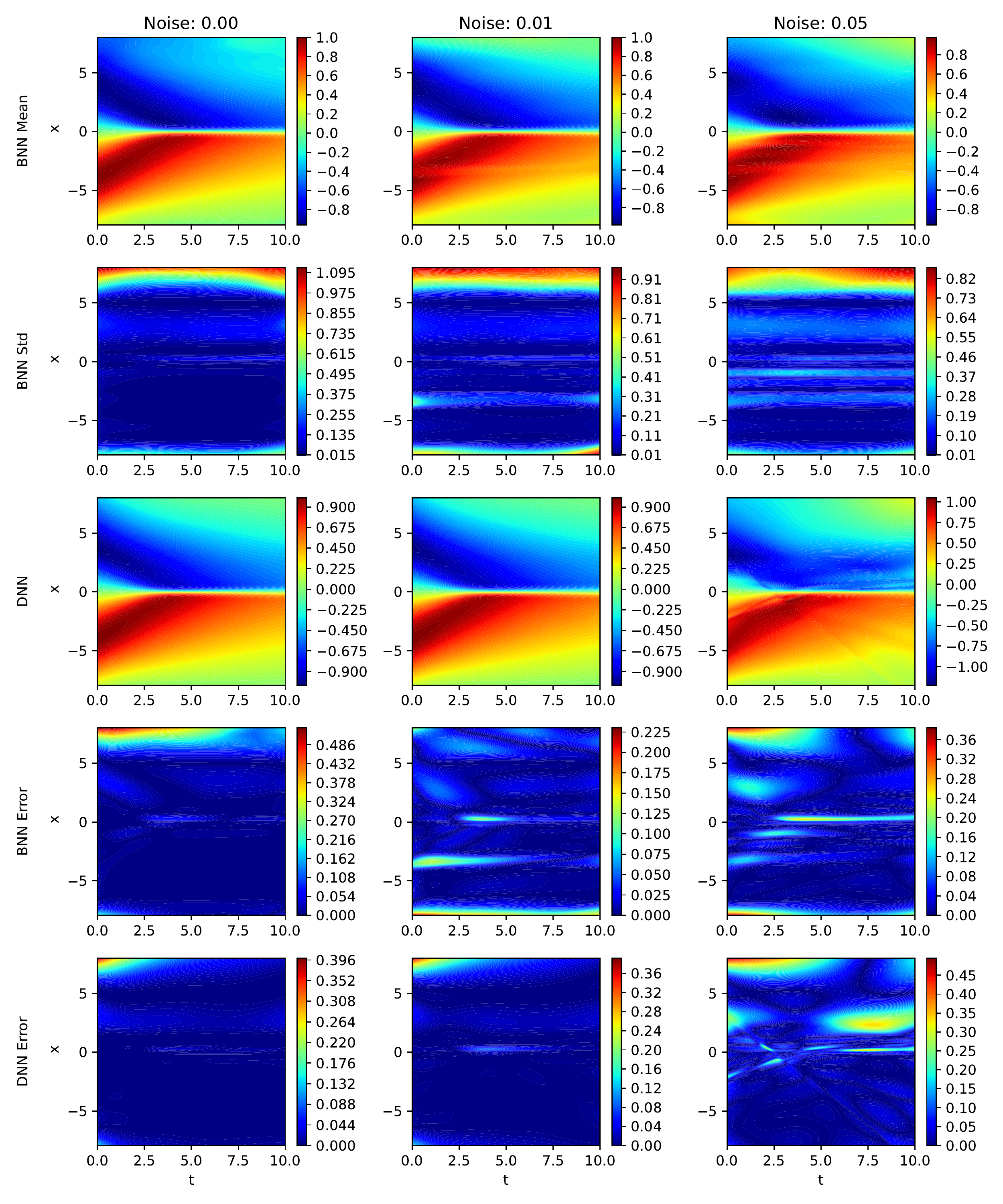}
  \caption{BNN predictive mean and standard deviation, DNN prediction, and absolute error between the BNN predictive mean, DNN and ground truth - Burgers Equation}
  \label{fig5}
\end{figure}

\begin{figure}[!h]
\hspace{-2.9cm}
  \includegraphics[width=1.5\textwidth]{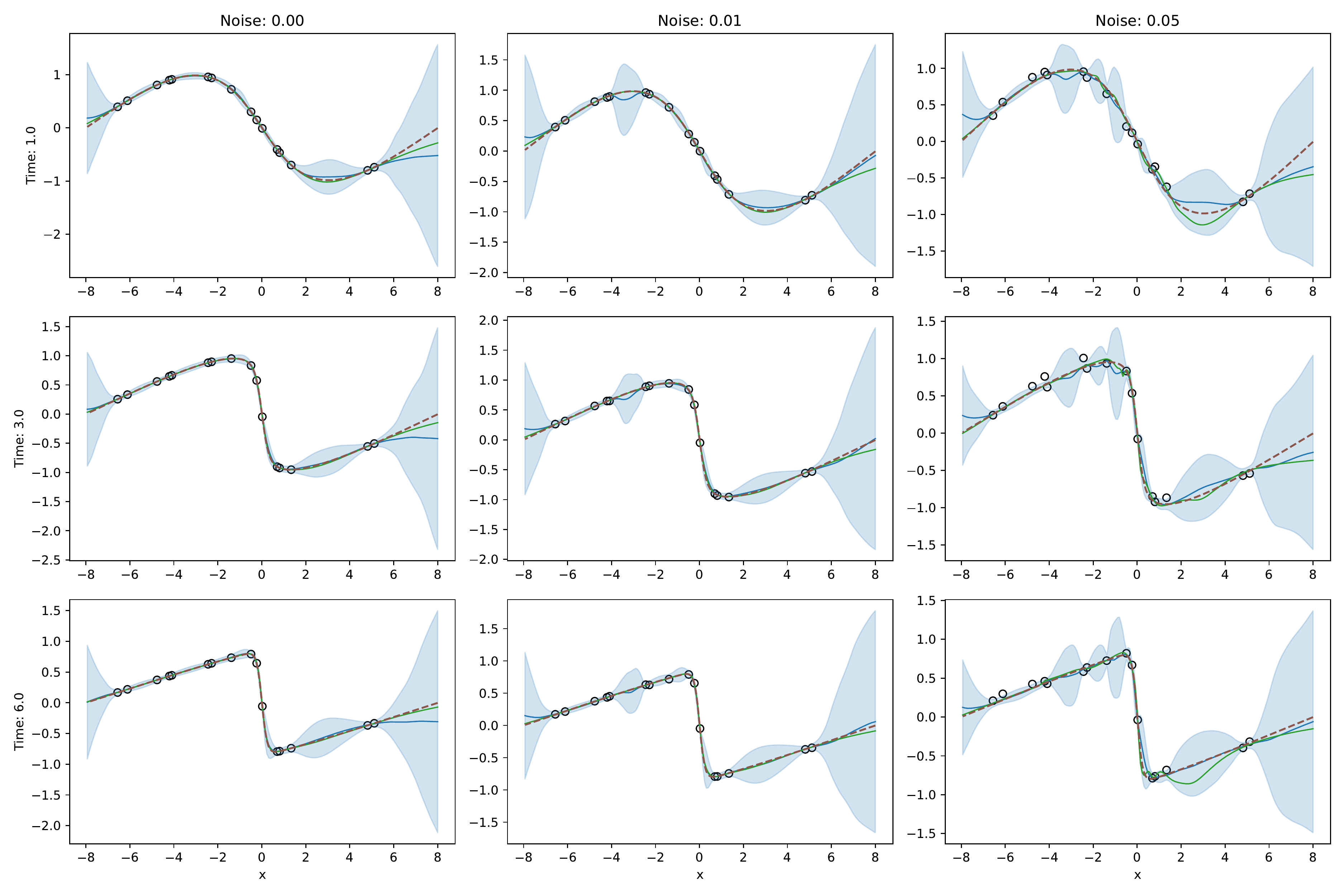}
  \caption{Predictions of the BNN (solid blue line $95\%$ confidence intervals) and DNN (solid green line) compared to the ground truth (dashed red line) at $t=1~s$, $3~s$ and $6~s$. The black circles represent sensor measurement data - Burgers Equation}
  \label{fig6}
\end{figure}

\noindent Let us consider Burgers equation example as originally proposed in \cite{JMLR:v19:18-046}:
\begin{equation}
    \frac{\partial u}{\partial t}=\lambda u\frac{\partial u}{\partial x}+\nu\frac{\partial^2 u}{\partial^2 x}
\end{equation}
\begin{equation}
\begin{cases}
    \,u(t=0,x)=-\sin(\pi x/8)\\[2pt]
    \,u(t,x=\pm8)=0
\end{cases}
\hspace{2cm}
\begin{cases}
    \,(t,x)\in\Omega=[0,10]\times[-8,+8]\\[2pt]
    \,\lambda=-1\hspace{1cm}\nu=0.1
\end{cases}
\end{equation}
Measurements are recorded every $\Delta t=0.2~s$ for each of the 16 sensors ($n=800$ training points), the ground truth solution is shown in figure \ref{fig4}. The predictions of the BNN and DNN are shown in figure \ref{fig5} and \ref{fig6}, and table \ref{table1} shows the RMSE between the BNN predictive mean, DNN prediction and ground truth over 1000 random input points. 

\FloatBarrier

\FloatBarrier

\begin{table}
\centering
\begin{tabular}{c|c|c}
Noise & BNN  & DNN \\\hline
$\epsilon=0$ & 0.0740 & \textbf{0.0283}\\\hline
$\epsilon\sim\mathcal{N}(0,0.01^2)$ & 0.0327 & \textbf{0.0298}\\\hline
$\epsilon\sim\mathcal{N}(0,0.05^2)$ & \textbf{0.0641} & 0.0919\\\hline
\end{tabular}
\caption{\centering\label{table1} RMSE of the Bayesian and Standard Deep Neural Network - Burgers Equation}
\end{table}

\noindent The dictionary of spatial derivative candidates, with the discovered coefficients using the STBLR trained on the BNN derivative data, is shown in table \ref{table2}. Table \ref{table3} shows the discovered coefficients obtained with the STOLS trained on the DNN derivative data (no uncertainty weighing). In the first two noise cases, although the DNN fits the measurement data better (the prediction RMSE on the test set is lower, as shown in table \ref{table1}), the BNN is able to better discover the underlying PDE, thanks to the well quantified uncertainty. For larger amount of noise ($\epsilon\sim\mathcal{N}(0,0.05^2)$), the STBLR trained on the BNN derivative data slightly underestimates the influence of both the non-linear term and the diffusion term in Burgers equation, and also slightly overestimates other derivative terms. However, it still outperforms the baseline comparison: the STOLS trained on the DNN data completely misses the diffusion term, and dramatically underestimate the non-linear term. 
\\\\
\noindent As shown in table \ref{table1}, the BNN is able to limit overfitting better than the standard DNN in cases with a larger amount of noise. In figure \ref{fig6}, we notice how the standard DNN (green line) becomes oscillatory around noisy datapoints, while the BNN captures the underlying pattern
correctly. The later also provides meaningful confidence intervals: As seen on figure \ref{fig5} and \ref{fig6}, the confidence intervals become wider (i.e. higher uncertainty) in space-time regions far from any measurement data, and the prediction errors are concentrated in areas of high uncertainty. Consequently, the BNN is able to produce better
derivative and uncertainty estimates, ultimately yielding higher PDE discovery accuracy than the DNN coupled with STOLS.

\FloatBarrier

\begin{figure}[!h]
\hspace{-2.9cm}
  \includegraphics[width=1.5\textwidth]{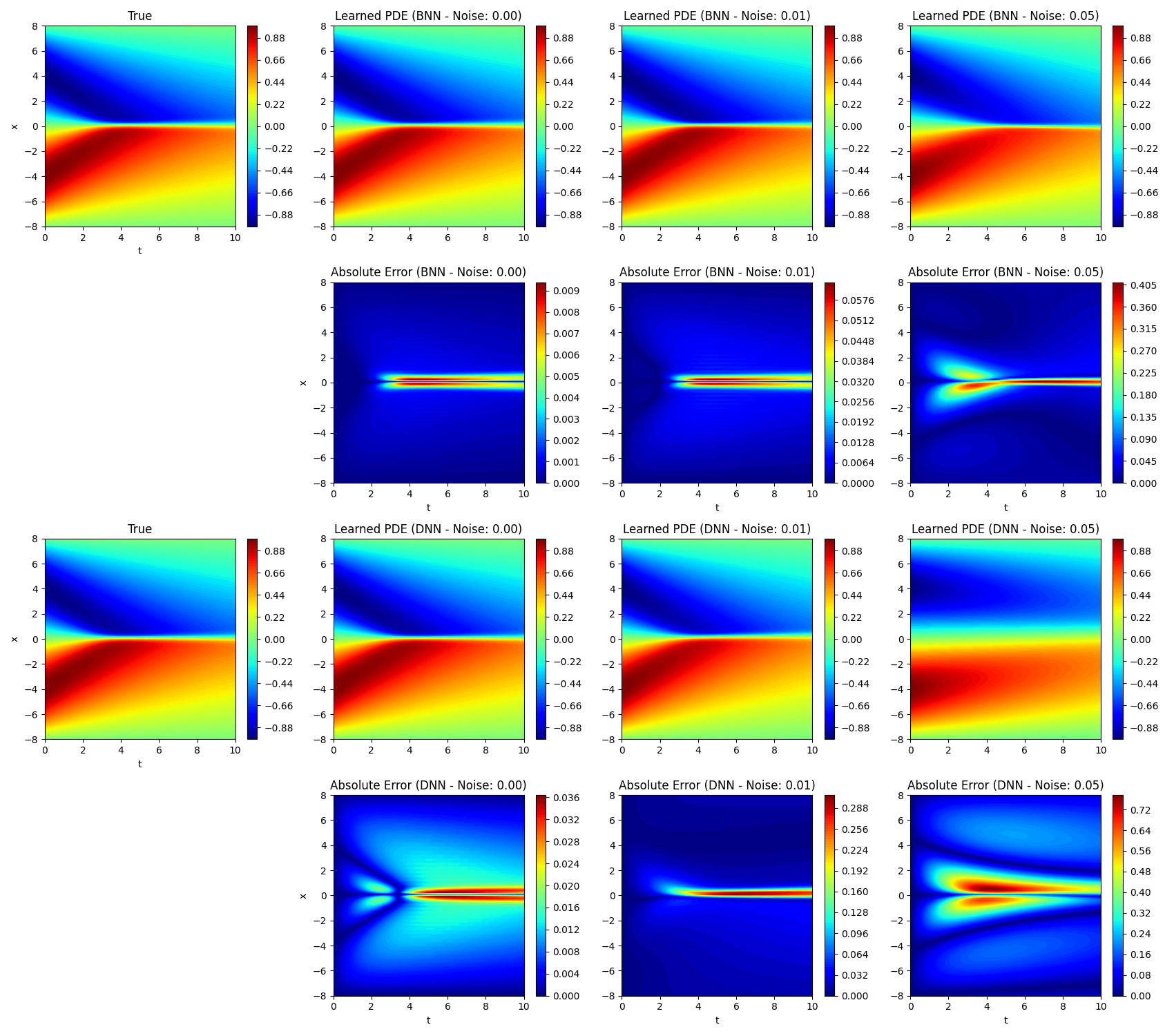}
  \caption{Learned solutions for the Burgers equations. First and third row show the solutions of the discovered PDE, with the BNN/STBLR and the DNN/STOLS, respectively. Second and fourth row show the absolute error with respect to the ground truth for the BNN/STBLR and the DNN/STOLS, respectively}
  \label{learntburgers}
\end{figure}

\noindent Figure \ref{learntburgers} shows the learned dynamics (i.e. the solution of the discovered PDE) for the STBLR trained on the BNN data, compared with the STOLS trained on the DNN data. For every noise cases, the BNN/STBLR does a remarkable job at learning the PDE dynamics accurately. The DNN/STOLS is able to capture the dynamics fairly accurately for cases with little noise, but it clearly fails to capture the shock that occurs over time around $x=0$ in the most noisy case. 

\begin{table}
\centering
\begin{tabular}{c|c|c|c|c}
Candidates & Ground Truth & Noiseless  & $\epsilon\sim\mathcal{N}(0,0.01^2)$ & $\epsilon\sim\mathcal{N}(0,0.05^2)$ \\\hline
 $u$ & $0.0$ & $0.0$ & $0.0$ & $-0.0415$ \\[-5pt] 
 $u_{x}$ & $0.0$ & $0.0$ & $0.0$ & $0.0165$ \\[-5pt] 
 $u_{xx}$ & $0.1$ & $0.0980$ & $0.0871$ & $0.0470$ \\[-5pt] 
 $u_{xxx}$ & $0.0$ & $0.0$ & $0.0$ & $0.0$ \\[-5pt] 
 $u_{xxxx}$ & $0.0$ & $0.0$ & $0.0$ & $0.0$ \\[-5pt] 
 $uu_{x}$ & $-1.0$ & $-0.9986$ & $-0.9851$ & $-0.6884$ \\[-5pt] 
 $uu_{xx}$ & $0.0$ & $0.0$ & $0.0$ & $-0.0253$ \\[-5pt] 
 $uu_{xxx}$ & $0.0$ & $0.0$ & $0.0$ & $0.0$ \\[-5pt] 
 $uu_{xxxx}$ & $0.0$ & $0.0$ & $0.0$ & $0.0$ \\[-5pt] 
 $u^2$ & $0.0$ & $0.0$ & $0.0$ & $-0.0202$ \\[-5pt] 
 $u^2_{x}$ & $0.0$ & $0.0$ & $0.0$ & $-0.0268$ \\\hline
$e_C$ ($\ell^2$ norm) & $0.0$ & $\mathbf{0.0024}$ & $\mathbf{0.0198}$ & $\mathbf{0.3220}$\\[-5pt] 
$e_L$ ($\ell^2$ norm) & $0.0$ & $\mathbf{0.1352}$ & $\mathbf{0.9446}$ & $\mathbf{6.3256}$
\end{tabular}
\caption{\centering\label{table2}Dictionary of candidate derivatives and discovered coefficients using STBLR and BNN derivative data ($\delta=0.005$) - Burgers Equation}
\end{table}

\begin{table}
\centering
\begin{tabular}{c|c|c|c|c}
Candidates & Ground Truth & Noiseless  & $\epsilon\sim\mathcal{N}(0,0.01^2)$ & $\epsilon\sim\mathcal{N}(0,0.05^2)$ \\\hline
 $u$ & $0.0$ & $0.0$ & $-0.0134$ & $-0.0609$ \\[-5pt] 
 $u_{x}$ & $0.0$ & $0.0$ & $-0.0156$ & $-0.0354$ \\[-5pt] 
 $u_{xx}$ & $0.1$ & $0.0918$ & $0.0921$ & $0.0$ \\[-5pt] 
 $u_{xxx}$ & $0.0$ & $0.0$ & $0.0$ & $0.0$ \\[-5pt] 
 $u_{xxxx}$ & $0.0$ & $0.0$ & $0.0$ & $0.0$ \\[-5pt] 
 $uu_{x}$ & $-1.0$ & $-0.9685$ & $-0.9181$ & $-0.2041$ \\[-5pt] 
 $uu_{xx}$ & $0.0$ & $0.0$ & $-0.0178$ & $0.0$ \\[-5pt] 
 $uu_{xxx}$ & $0.0$ & $0.0$ & $0.0$ & $0.0$ \\[-5pt] 
 $uu_{xxxx}$ & $0.0$ & $0.0$ & $0.0$ & $0.0$ \\[-5pt] 
 $u^2$ & $0.0$ & $0.0$ & $0.0$ & $0.0$ \\[-5pt] 
 $u^2_{x}$ & $0.0$ & $0.0$ & $-0.0265$ & $0.0$ \\\hline
$e_C$ ($\ell^2$ norm) & $0.0$ & $0.0326$ & $0.0906$ & $0.8052$\\[-5pt] 
$e_L$ ($\ell^2$ norm) & $0.0$ & $0.9962$ & $4.0730$ & $20.4654$
\end{tabular}
\caption{\centering\label{table3} Dictionary of candidate derivatives and discovered coefficients using STOLS and DNN derivative data ($\delta=0.005$) - Burgers Equation}
\end{table}

\FloatBarrier

\subsection{\label{kdv} Korteweg-De Vries Equation}

\begin{figure}[!h]
\centering
  \centering
  \includegraphics[width=1\textwidth]{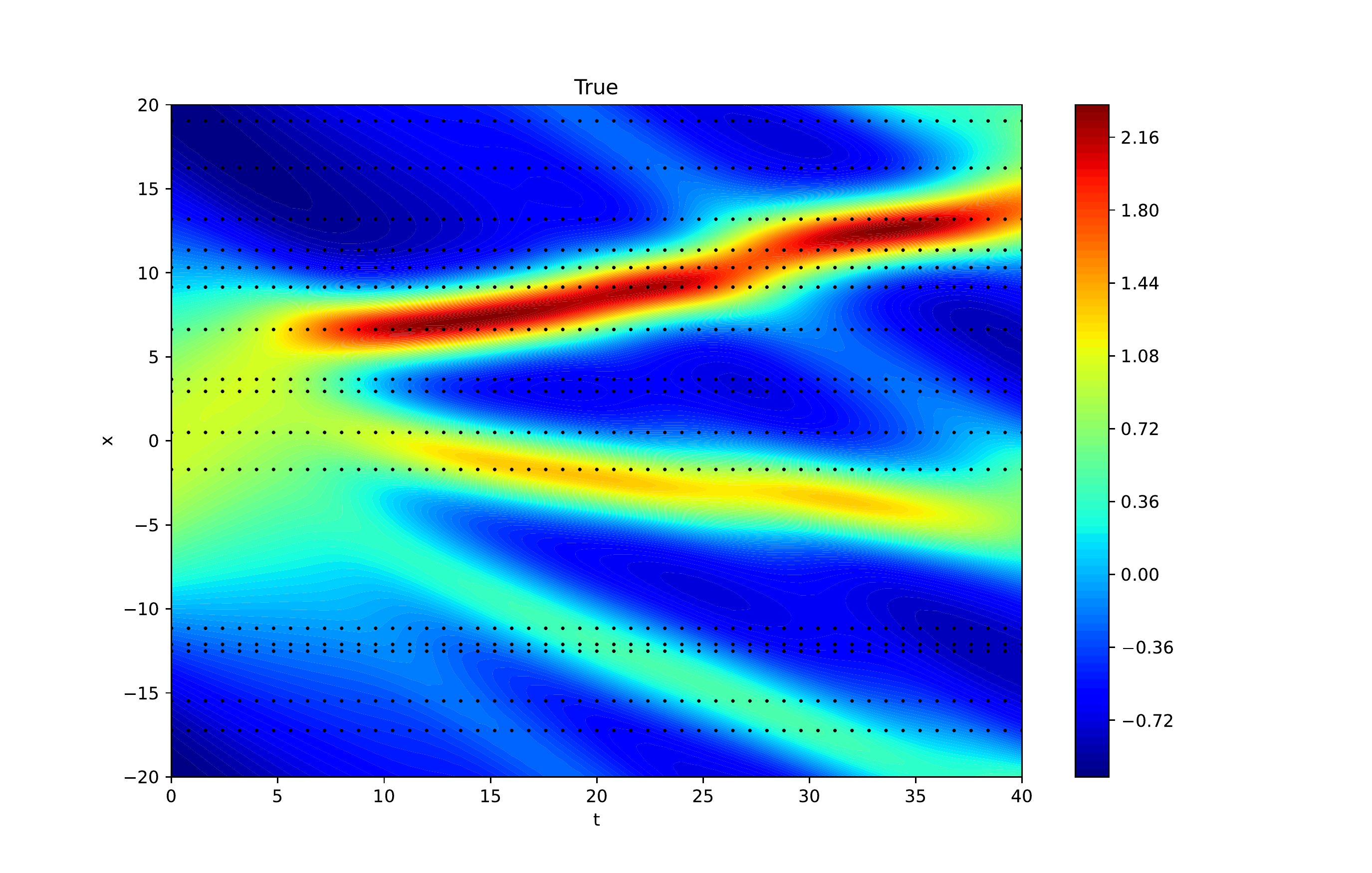}
  \caption{Ground truth solution for KdV equation (black dots represent measured time-space points)}
  \label{fig7}
\end{figure}

\begin{figure}[!h]
\vspace{-2cm}
\hspace{-2.8cm}
  \includegraphics[width=1.4\textwidth]{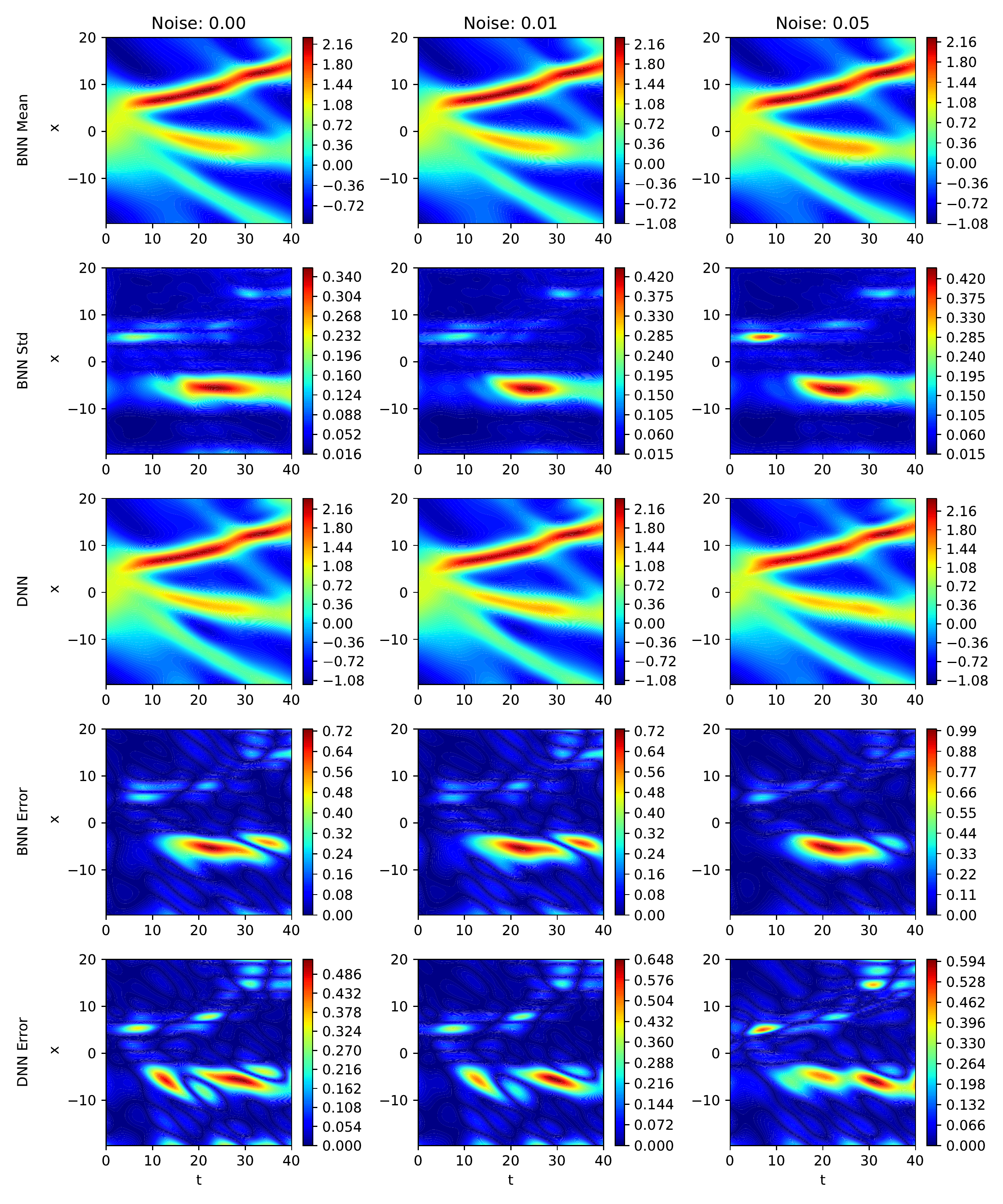}
  \caption{BNN predictive mean and standard deviation, DNN prediction, and absolute error between the BNN predictive mean, DNN and ground truth - KdV Equation}
  \label{fig8}
\end{figure}

\begin{figure}[!h]
\hspace{-2.9cm}
  \includegraphics[width=1.5\textwidth]{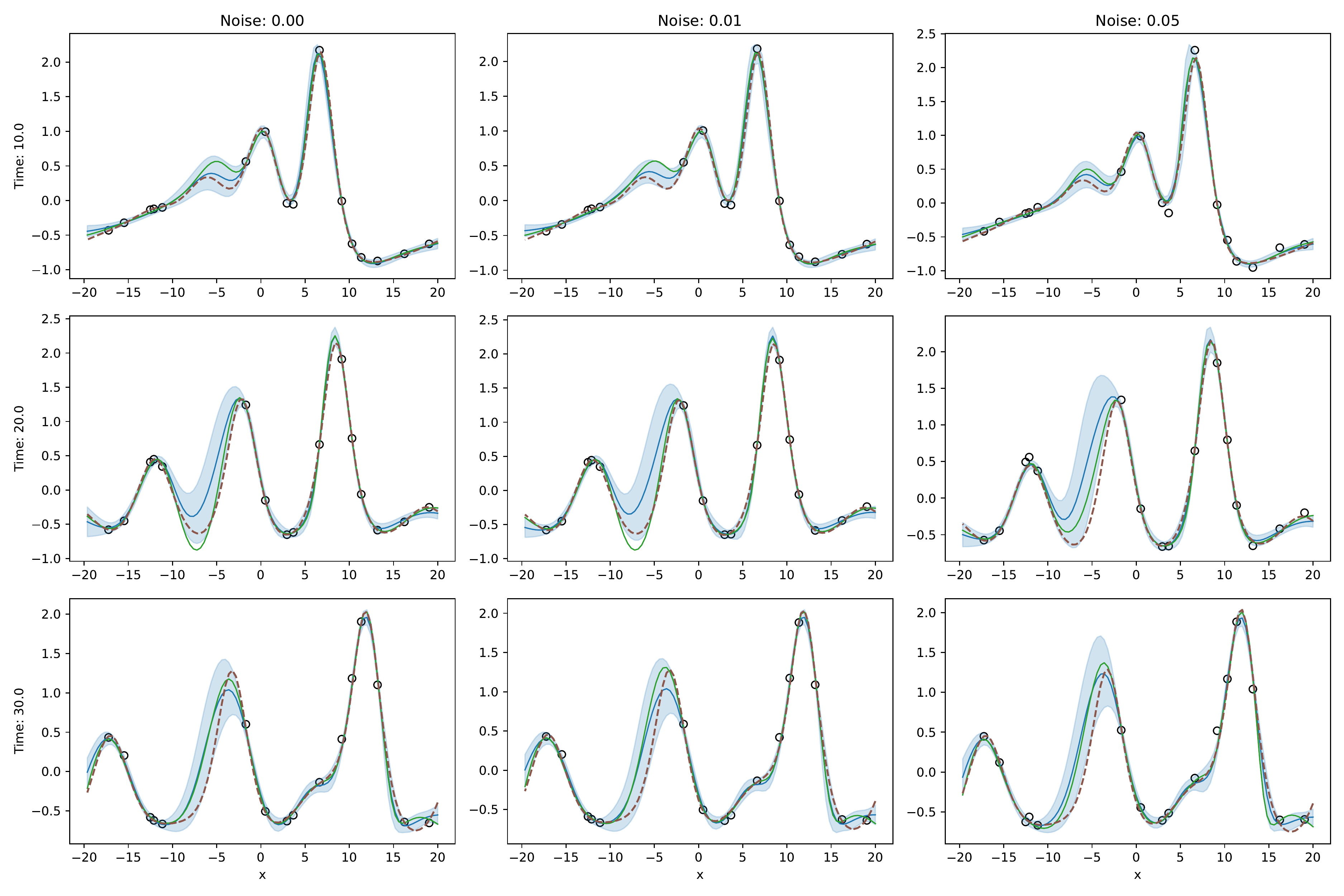}
  \caption{Predictions of the BNN (solid blue line $95\%$ confidence intervals) and DNN (solid green line) compared to the ground truth (dashed red line) at $t=10~s$, $20~s$ and $30~s$. The black circles represent sensor measurement data - KdV Equation}
  \label{fig9}
\end{figure}

\noindent We now consider the Korteweg-De Vries (KdV) equation example as originally proposed in \cite{JMLR:v19:18-046}:
\begin{equation}
    \frac{\partial u}{\partial t}=\lambda u\frac{\partial u}{\partial x}+\beta\frac{\partial^3u}{\partial^3 x}
\end{equation}
\begin{equation}
\begin{cases}
    \,u(t=0,x)=\cos(-\pi x/20)\\[2pt]
    \,(t,x)\in\Omega=[0,40]\times[-20,+20]
\end{cases}
\hspace{2cm}
\begin{cases}
    \,\lambda=-1\\[2pt]
    \,\beta=-1
\end{cases}
\end{equation}
Measurements are recorded every $\Delta t=0.8~s$ for each of the 16 sensors ($n=800$ training points), the ground truth solution is shown in figure \ref{fig7}. The predictions of the BNN and DNN are shown in figure \ref{fig8} and \ref{fig9}.
\\\\

\FloatBarrier

\noindent Table \ref{table4} shows the prediction RMSE for both the BNN and DNN. The DNN predictions outperform the BNN, but then fails to recover the ground truth PDE, as shown in table \ref{table5} and \ref{table6}. Even in the cases with limited to no noise at all, the DNN does a poor job. While it correctly identifies the third order spatial derivative and the non-linear term, it dramatically underestimate their influence. Conversely, the STBLR trained on the BNN derivative data not only identifies the derivatives correctly, but also predicts their coefficients much more accurately. This example shows that the accuracy of the BNN predictions is independent from its capability to discover the governing PDE. Indeed, the BNN doesn't have to be accurate everywhere, it only takes a few correct derivative estimations, with well quantified uncertainty, to find the ground truth PDE.
\\\\

\begin{table}
\centering
\begin{tabular}{c|c|c}
Noise & BNN  & DNN \\\hline
$\epsilon=0$ & 0.1388 & \textbf{0.1018}\\\hline
$\epsilon\sim\mathcal{N}(0,0.01^2)$ & 0.1437 & \textbf{0.1068}\\\hline
$\epsilon\sim\mathcal{N}(0,0.05^2)$ & 0.1724 & \textbf{0.1158}\\\hline
\end{tabular}
\caption{\centering\label{table4}RMSE of the Bayesian and Standard Deep Neural Network - KdV Equation}
\end{table}

\noindent In the case with large amount of noise ($\epsilon\sim\mathcal{N}(0,0.05^2)$), the STBLR trained on the BNN derivative data still discovers the third order spatial derivative, and the non-linear term, with satisfactory accuracy; clearly outperforming the DNN. Indeed, the STOLS trained on the DNN derivative data misses the third order spatial derivative and dramatically underestimates the non-linear term.

\begin{table}
\centering
\begin{tabular}{c|c|c|c|c}
Candidates & Ground Truth & Noiseless  & $\epsilon\sim\mathcal{N}(0,0.01^2)$ & $\epsilon\sim\mathcal{N}(0,0.05^2)$ \\\hline
 $u$ & $0$ & $0.0$ & $0.0$ & $0.0$ \\[-5pt] 
 $u_{x}$ & $0$ & $0.0$ & $0.0$ & $0.0$ \\[-5pt] 
 $u_{xx}$ & $0$ & $0.0$ & $0.0$ & $0.0$ \\[-5pt] 
 $u_{xxx}$ & $-1$ & $-0.9296$ & $-0.9254$ & $-0.8653$ \\[-5pt] 
 $u_{xxxx}$ & $0$ & $0.0$ & $0.0$ & $0.0$ \\[-5pt] 
 $uu_{x}$ & $-1$ & $-0.8199$ & $-0.8160$ & $-0.8303$ \\[-5pt] 
 $uu_{xx}$ & $0$ & $0.0$ & $0.0$ & $0.0$ \\[-5pt] 
 $uu_{xxx}$ & $0$ & $0.0$ & $0.0$ & $0.0$ \\[-5pt] 
 $uu_{xxxx}$ & $0$ & $0.0$ & $0.0$ & $0.0$ \\[-5pt] 
 $u^2$ & $0$ & $0.0$ & $0.0$ & $0.0$ \\[-5pt] 
 $u^2_{x}$ & $0$ & $0.0$ & $0.0$ & $0.0$ \\\hline
$e_C$ ($\ell^2$ norm) & $0.0$ & $\mathbf{0.1934}$ & $\mathbf{0.1986}$ & $\mathbf{0.2166}$\\[-5pt] 
$e_L$ ($\ell^2$ norm) & $0.0$ & $\mathbf{119.2665}$ & $\mathbf{120.8644}$ & $\mathbf{108.1781}$
\end{tabular}
\caption{\centering\label{table5}Dictionary of candidate derivatives and discovered coefficients using STBLR and BNN derivative data ($\delta=0.05$) - KdV Equation}
\end{table}

\begin{table}
\centering
\begin{tabular}{c|c|c|c|c}
Candidates & Ground Truth & Noiseless  & $\epsilon\sim\mathcal{N}(0,0.01^2)$ & $\epsilon\sim\mathcal{N}(0,0.05^2)$ \\\hline
 $u$ & $0$ & $0.0$ & $0.0$ & $0.0$ \\[-5pt] 
 $u_{x}$ & $0$ & $0.0$ & $0.0$ & $0.0$ \\[-5pt] 
 $u_{xx}$ & $0$ & $0.0$ & $0.0$ & $0.0$ \\[-5pt] 
 $u_{xxx}$ & $-1$ & $-0.2056$ & $-0.2061$ & $0.0$ \\[-5pt] 
 $u_{xxxx}$ & $0$ & $0.0$ & $0.0$ & $0.0$ \\[-5pt] 
 $uu_{x}$ & $-1$ & $-0.3008$ & $-0.2961$ & $-0.1330$ \\[-5pt] 
 $uu_{xx}$ & $0$ & $0.0$ & $0.0$ & $0.0$ \\[-5pt] 
 $uu_{xxx}$ & $0$ & $0.0$ & $0.0$ & $0.0$ \\[-5pt] 
 $uu_{xxxx}$ & $0$ & $0.0$ & $0.0$ & $0.0$ \\[-5pt] 
 $u^2$ & $0$ & $0.0$ & $0.0$ & $0.0$ \\[-5pt] 
 $u^2_{x}$ & $0$ & $0.0$ & $0.0$ & $0.0$ \\\hline
$e_C$ ($\ell^2$ norm) & $0.0$ & $1.0583$ & $1.061$ & $1.3235$\\[-5pt] 
$e_L$ ($\ell^2$ norm) & $0.0$ & $233.8474$ & $234.1718$ & $232.9447$\\
\end{tabular}
\caption{\centering\label{table6} Dictionary of candidate derivatives and discovered coefficients using STOLS and DNN derivative data ($\delta=0.05$) - KdV Equation}
\end{table}

\FloatBarrier

\begin{figure}[!h]
\hspace{-2.9cm}
  \includegraphics[width=1.5\textwidth]{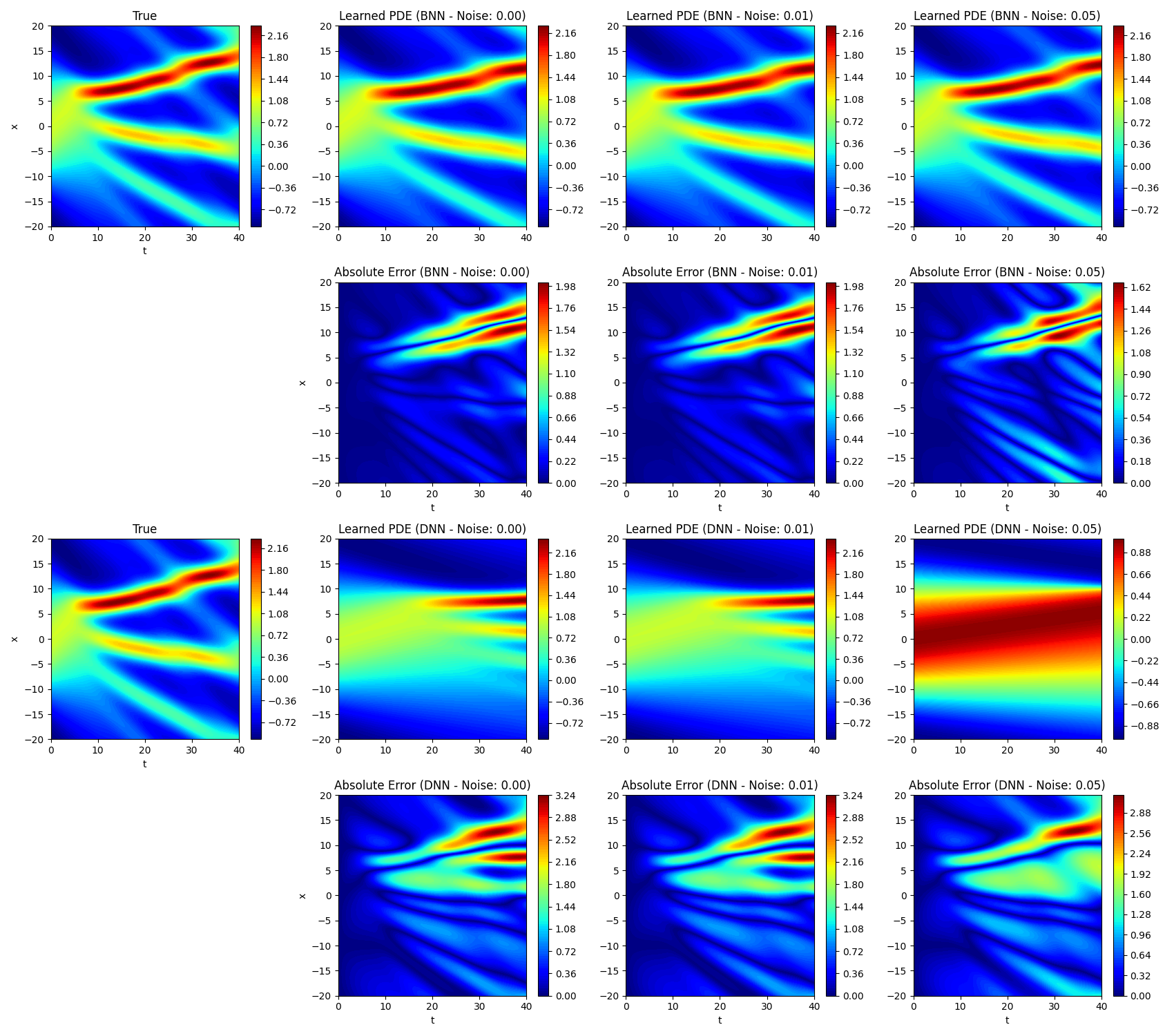}
  \caption{Learned solutions for the KdV equations. First and third row show the solutions of the discovered PDE, with the BNN/STBLR and the DNN/STOLS, respectively. Second and fourth row show the absolute error with respect to the ground truth for the BNN/STBLR and the DNN/STOLS, respectively}
  \label{learntKdV}
\end{figure}

\noindent Figure \ref{learntKdV} shows the learned dynamics for the STBLR trained on the BNN data, compared with the STOLS trained on the DNN data. For every noise case, the BNN/STBLR is able to learn the PDE dynamics with high accuracy (despite some errors localized around the wave front). Conversely, in each cases (and particularly the most noisy case), The DNN/STOLS dramatically fails to capture accurately the dynamics.

\FloatBarrier

\subsection{\label{heat} Heat Equation}
\begin{figure}[!h]
\centering
  \centering
  \includegraphics[width=1\textwidth]{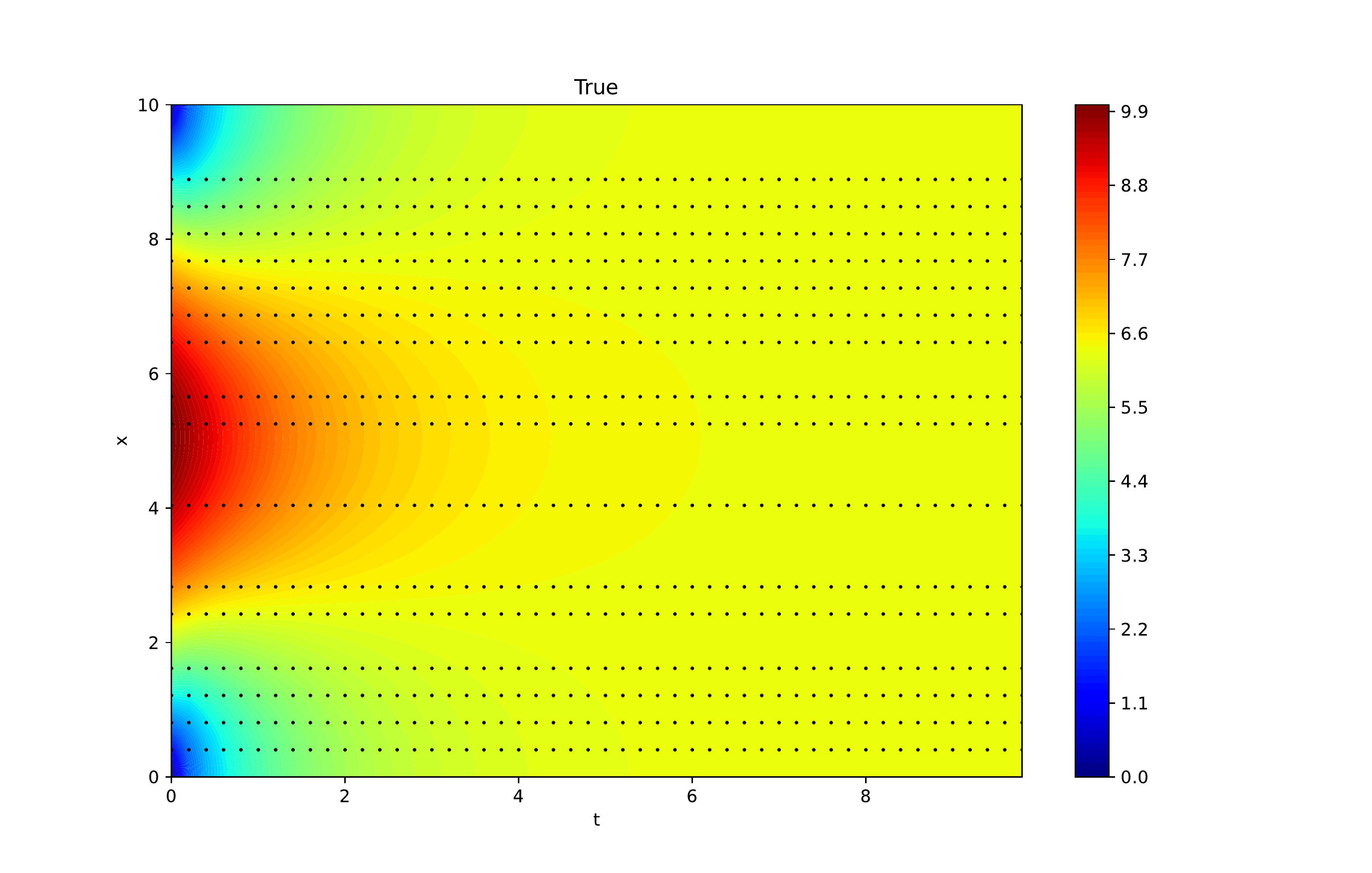}
  \caption{Ground truth solution for the heat equation (black dots represent measured time-space points)}
  \label{fig1}
\end{figure}

\begin{figure}[!h]
\vspace{-2cm}
\hspace{-2.8cm}
  \includegraphics[width=1.4\textwidth]{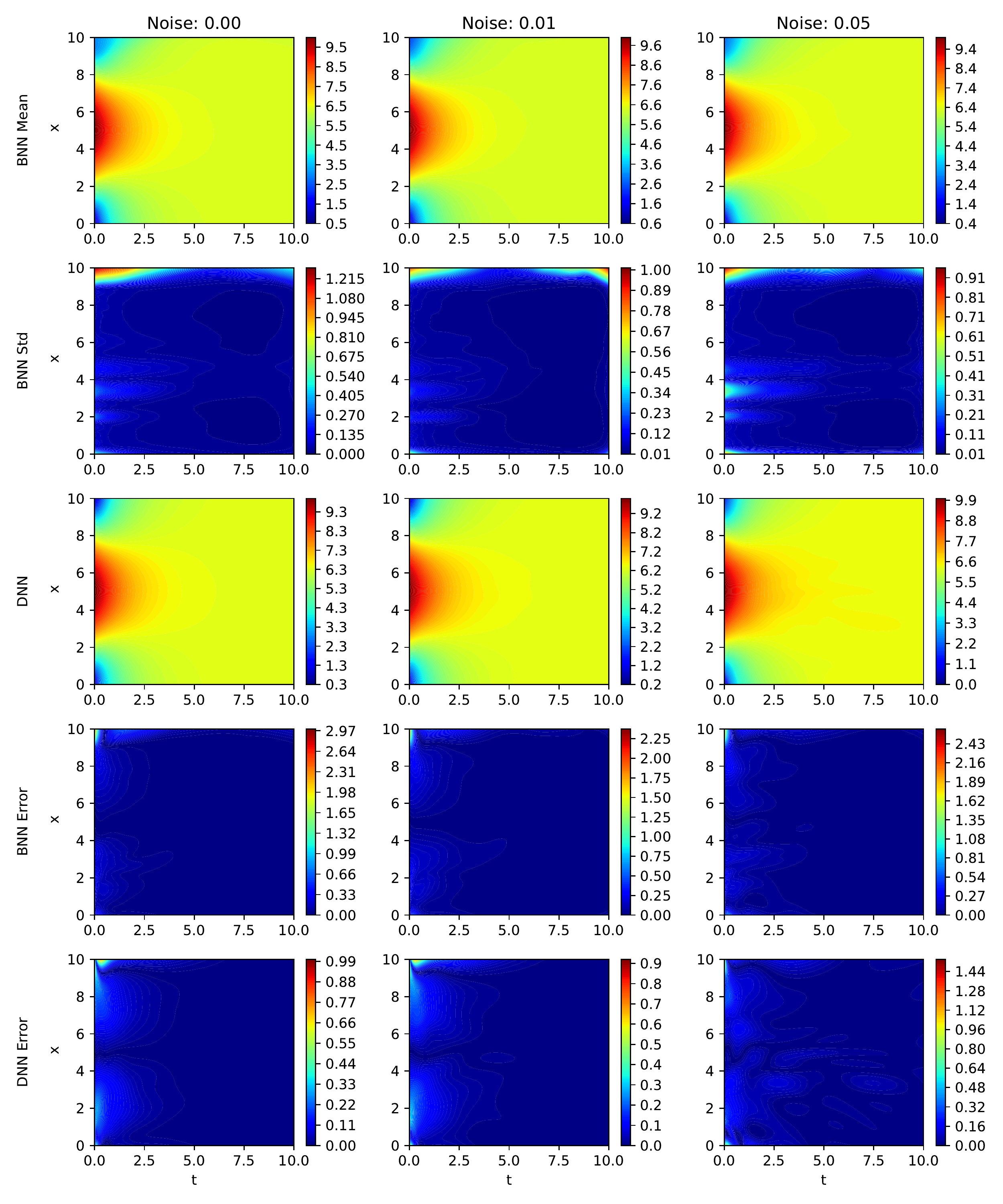}
  \caption{BNN predictive mean and standard deviation, DNN prediction, and absolute error between the BNN predictive mean, DNN and ground truth - Heat equation}
  \label{fig2}
\end{figure}

\begin{figure}
\hspace{-2.9cm}
  \includegraphics[width=1.5\textwidth]{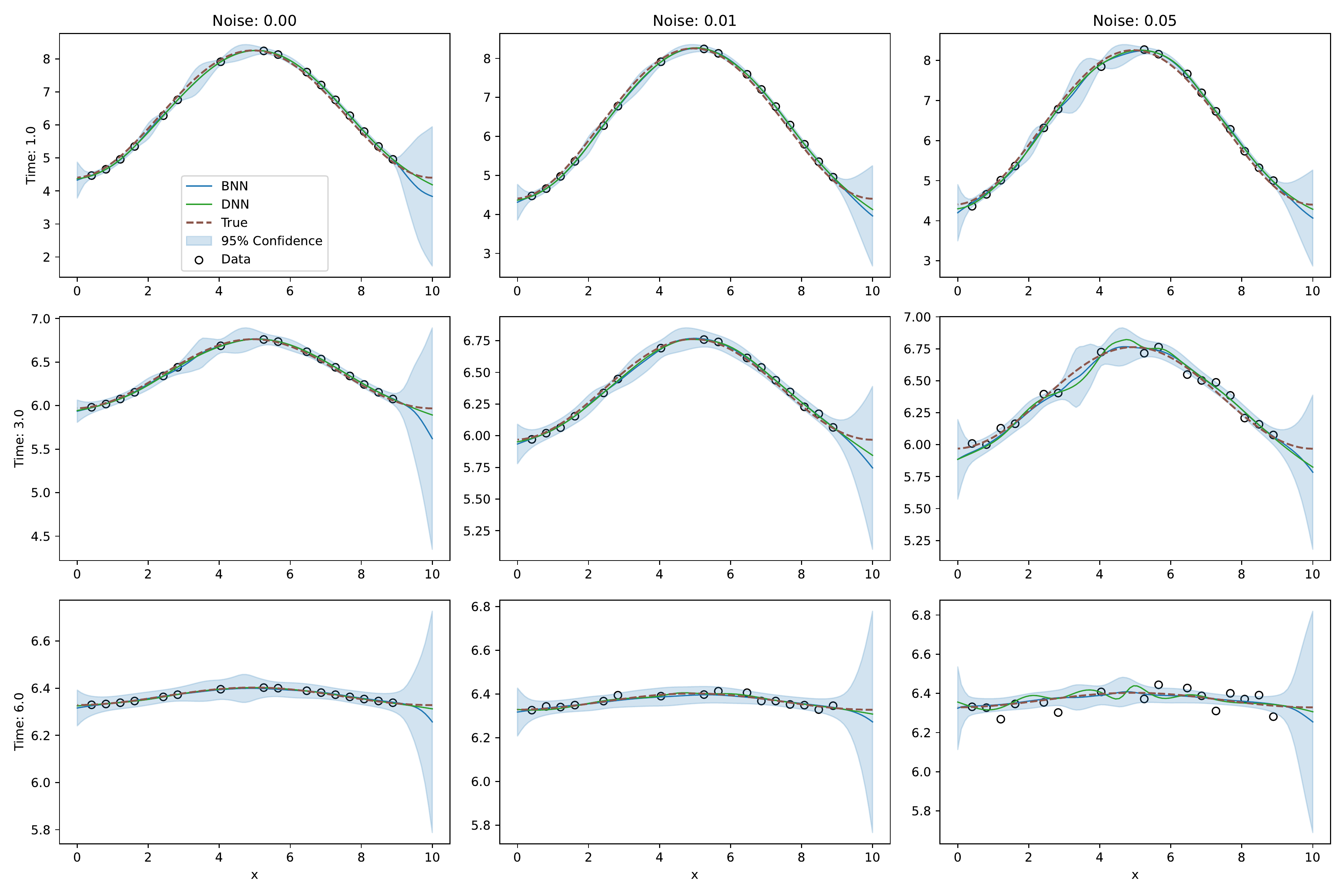}
  \caption{Predictions of the BNN (solid blue line $95\%$ confidence intervals) and DNN (solid green line) compared to the ground truth (dashed red line) at $t=1~s$, $3~s$ and $6~s$. The black circles represent sensor measurement data - Heat equation}
  \label{fig3}
\end{figure}

\noindent Finally, we consider the following 1D heat equation:
\begin{equation}
    \frac{\partial u}{\partial t}=\nu\frac{\partial^2 u}{\partial^2 x}
\end{equation}
\begin{equation}
\begin{cases}
    \,u(t=0,x)=10\cos(\pi(x-5)/10)\\[2pt]
    \,\displaystyle\frac{\partial u}{\partial x}(t,x=0,10)=0
\end{cases}
\hspace{1cm}
\begin{cases}
    \,(t,x)\in\Omega=[0,10]\times[0,10]\\[2pt]
    \,\nu=2
\end{cases}
\end{equation}
The temperature is recorded every $\Delta t=0.2~s$ for each of the 16 sensors ($n=800$ training points), and the ground truth solution is shown in figure \ref{fig1}. The predictions of the BNN and DNN are shown in figure \ref{fig2} and \ref{fig3}.
\\\\

\FloatBarrier

\FloatBarrier

\FloatBarrier

\FloatBarrier

\begin{figure}[!h]
\hspace{-2.9cm}
  \includegraphics[width=1.5\textwidth]{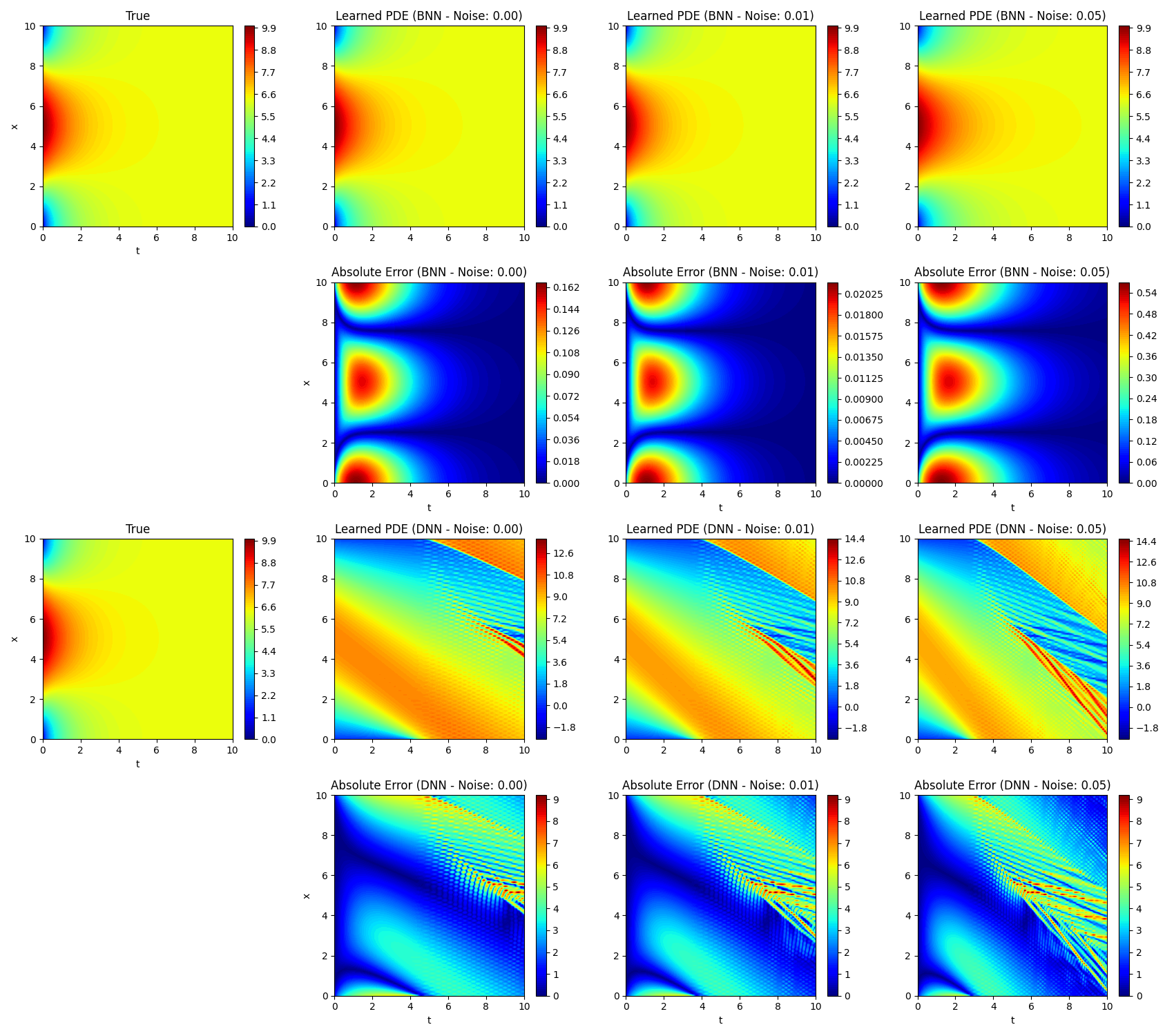}
  \caption{Learned solutions for the heat equation. First and third row show the solutions of the discovered PDE, with the BNN/STBLR and the DNN/STOLS, respectively. Second and fourth row show the absolute error with respect to the ground truth for the BNN/STBLR and the DNN/STOLS, respectively}
  \label{learntheat}
\end{figure}

\noindent Table \ref{table7} shows the prediction RMSE for both the BNN and DNN and figure \ref{fig3} represents the BNN and DNN predictions at different times. Similarly to previous examples, the standard DNN yields lower errors on the test set. However, as shown in table \ref{table8} and \ref{table9}, the STBLR trained on the BNN derivative data is able to recover the ground truth diffusion equation much better than the DNN/STOLS model. It accurately detects the second order spatial derivative (though with a slight underestimation of the diffusivity coefficient), while the DNN/STOLS models fails to find the diffusion term, and instead surprisingly discovers a squared first order derivative term that is not remotely comprised in the ground truth heat equation.
\\\\
\noindent Figure \ref{learntheat} shows the dynamics of the discovered PDE. Since the STBLR trained on the BNN data correctly identified the heat equation, the dynamics for each noise case is very consistent with the ground truth, although the heat diffusion is a little slower due to the underestimated diffusivity coefficient. For comparison, the non-linear PDE discovered using the baseline DNN/STOLS model leads to dynamics completely different from the ground truth. Note that in this later case, the numerical solutions of the discovered PDE are likely becoming unstable for $t>3$, but the dynamic is already wrong before this.

\begin{table}
\centering
\begin{tabular}{c|c|c}
Noise & BNN  & DNN \\\hline
$\epsilon=0$ & 0.0593 & \textbf{0.0203}\\\hline
$\epsilon\sim\mathcal{N}(0,0.01^2)$ & 0.0395 & \textbf{0.0243}\\\hline
$\epsilon\sim\mathcal{N}(0,0.05^2)$ & 0.0399 & \textbf{0.0310}\\\hline

\end{tabular}
\caption{\centering\label{table7}RMSE in the predicted heat equation response of the Bayesian and Standard Deep Neural Networks}
\end{table}

\begin{table}
\centering
\begin{tabular}{c|c|c|c|c}
Candidates & Ground Truth & Noiseless  & $\epsilon\sim\mathcal{N}(0,0.01^2)$ & $\epsilon\sim\mathcal{N}(0,0.05^2)$ \\\hline
 $u$ & $0$ & $0.0$ & $0.0$ & $0.0$ \\[-5pt] 
 $u_{x}$ & $0$ & $0.0$ & $0.0$ & $0.0$ \\[-5pt] 
 $u_{xx}$ & $2$ & $1.8089$ & $1.9742$ & $1.4124$ \\[-5pt] 
 $u_{xxx}$ & $0$ & $0.0$ & $0.0$ & $0.0$ \\[-5pt] 
 $u_{xxxx}$ & $0$ & $0.0$ & $0.0$ & $0.0$ \\[-5pt] 
 $uu_{x}$ & $0$ & $0.0$ & $0.0$ & $0.0$ \\[-5pt] 
 $uu_{xx}$ & $0$ & $0.0$ & $0.0$ & $0.0$ \\[-5pt] 
 $uu_{xxx}$ & $0$ & $0.0$ & $0.0$ & $0.0$ \\[-5pt] 
 $uu_{xxxx}$ & $0$ & $0.0$ & $0.0$ & $0.0$ \\[-5pt] 
 $u^2$ & $0$ & $0.0$ & $0.0$ & $0.0$ \\[-5pt] 
 $u^2_{x}$ & $0$ & $0.0$ & $0.0$ & $0.0$ \\\hline
$e_C$ ($\ell^2$ norm) & $0.0$ & $\mathbf{0.1911}$ & $\mathbf{0.0258}$ & $\mathbf{0.5876}$\\[-5pt] 
$e_L$ ($\ell^2$ norm) & $0.0$ & $\mathbf{7.8201}$ & $\mathbf{0.9894}$ & $\mathbf{28.7470}$\\
\end{tabular}
\caption{\centering\label{table8}Dictionary of candidate derivatives and discovered coefficients using STBLR and BNN derivative data ($\delta=0.02$) - Heat Equation}
\end{table}

\begin{table}
\centering
\begin{tabular}{c|c|c|c|c}
Candidates & Ground Truth & Noiseless  & $\epsilon\sim\mathcal{N}(0,0.01^2)$ & $\epsilon\sim\mathcal{N}(0,0.05^2)$ \\\hline
 $u$ & $0$ & $0.0$ & $0.0$ & $0.0$ \\[-5pt] 
 $u_{x}$ & $0$ & $0.0$ & $0.0$ & $0.0$ \\[-5pt] 
 $u_{xx}$ & $2$ & $0.0$ & $0.0$ & $0.0$ \\[-5pt] 
 $u_{xxx}$ & $0$ & $0.0$ & $0.0$ & $0.0$ \\[-5pt] 
 $u_{xxxx}$ & $0$ & $0.0$ & $0.0$ & $0.0$ \\[-5pt] 
 $uu_{x}$ & $0$ & $0.0$ & $0.0$ & $0.0$ \\[-5pt] 
 $uu_{xx}$ & $0$ & $0.0$ & $0.0$ & $0.0$ \\[-5pt] 
 $uu_{xxx}$ & $0$ & $0.0$ & $0.0$ & $0.0$ \\[-5pt] 
 $uu_{xxxx}$ & $0$ & $0.0$ & $0.0$ & $0.0$ \\[-5pt] 
 $u^2$ & $0$ & $0.0$ & $0.0$ & $0.0$ \\[-5pt] 
 $u^2_{x}$ & $0$ & $0.0376$ & $0.0461$ & $0.0601$ \\\hline
$e_C$ ($\ell^2$ norm) & $0.0$ & $2.0004$ & $2.0005$ & $2.0009$\\[-5pt] 
$e_L$ ($\ell^2$ norm) & $0.0$ & $401.6347$ & $404.4451$ & $409.0607$
\end{tabular}
\caption{\centering\label{table9} Dictionary of candidate derivatives and discovered coefficients using STOLS and DNN derivative data ($\delta=0.02$) - Heat Equation}
\end{table}

\FloatBarrier

\section{Conclusion}

\noindent Over the three application examples, the sequential threshold Bayesian linear regression model, used in conjunction with the Bayesian neural network expected derivatives and variance, clearly provides much more accurate results than the frequentist baseline comparison model (DNN/STOLS combination). Our proposed BNN/STBLR model is able to accurately discover both linear and non-linear PDEs with excellent accuracy; even in the noisiest cases and when assuming a large set of candidate derivatives. In noisy cases, the BNN is able to limit overfitting, without having to rely on data-demanding validation methods and specific tuning of the network hyperparameters. The BNN also provides valuable uncertainty quantification, which allows for discarding inaccurate derivative estimations; thus furnishing much better PDE discovery performance than its frequentist counterpart.
\\\\
\noindent While relying on BNNs to approximate the physical quantities of interest limits overfitting and helps quantifying uncertainty over the set of candidate derivatives, a potential caveat of this method is the need for approximate inference when training the neural network. First, Monte-Carlo methods, here used for sampling the posterior, are notoriously inefficient, and may be expensive when the data are more abundant (although in this case, fine-tuning a standard deep neural network may be feasible and preferable). Secondly, generating the dataset of candidate derivatives is computationally intensive: for each set of weights sampled from the posterior,  the weights are loaded into the neural network, a forward pass is performed, and the network is differentiated multiple times with respect to the inputs. Then, the derivatives are averaged over all the posterior samples. This process may be time consuming especially if a large number of weight samples is necessary. Approximating the posterior with a simpler (tractable) distribution along with variational inference may alleviate this burden, but also results in less informative derivative uncertainty quantification, and thus a less accurate PDE coefficient discovery.
\\\\
In this paper, we presented a framework for PDE discovery fully based on Bayesian inference, by combining Bayesian neural network with Bayesian linear regression. The use of BNNs for interpolating sparse physical measurements allows for accurately approximating the full physical system response, with well quantified uncertainty in regions where the data is seldom and noisy; ultimately allowing for more accurate PDE discovery. We believe that Bayesian methods can play a significant role in the field of PDE and dynamical system discovery, and this paper is a new contribution towards such a direction. 

\section{Acknowledgements}

\noindent Christopher J. Earls was supported by the Army Research Office (ARO) Biomathematics program, grant W911NF-18-1-0351.

\clearpage
\bibliography{references}

\end{document}